\documentclass{statsoc}

\RequirePackage[OT1]{fontenc}
\RequirePackage[colorlinks]{hyperref}
\RequirePackage{hypernat}

\usepackage{natbib}

\usepackage{latexsym,enumerate}
\usepackage{amsmath,amsopn, amstext, amscd, amsfonts, amssymb, euler,pdfsync, amsbsy }
\usepackage{mathrsfs}
\usepackage{fullpage}
\usepackage{graphicx}

\newtheorem{proposition}{Proposition}
\newtheorem{lemma}{Lemma}
\newtheorem{theorem}{Theorem}

\newtheorem{corollary}{Corollary}

\let\epsilon=\varepsilon
\let\eps=\epsilon
\let\phi=\varphi

\let\tilde=\widetilde

\newcommand{\field}[1]{\mathbb{#1}}
\newcommand{\R}{\field{R}}

\newcommand{\M}{{\mathcal M}}
\newcommand{\D}{\mathcal{D}}
\newcommand{\cI}{\mathcal{I}}

\newcommand{\cL}{\mathcal{L}}

\newcommand{\cT}{\mathcal{T}}
\newcommand{\cB}{{\mathcal B}}

\newcommand{\beqn}{\begin{equation}}
\newcommand{\eeqn}{\end{equation}}
\newcommand\eref[1]{(\ref{#1})}

\def\E{{\field E}}
\def\I{{\field I}}
\def\P{{\field P}}

\def\lan{\lambda_n}

\def\hf{\hat f}

\title[LOL]{Learning Out of Leaders}
\author[Mougeot {\it et al.}]{Mathilde Mougeot}
\address{Universit\'e Paris-Diderot, CNRS LPMA, 175 rue du Chevaleret, 75013
Paris, France.} \email{mathilde.mougeot@univ-paris-diderot.fr}
\author{Dominique Picard}
\address{Universit\'e Paris-Diderot, CNRS LPMA, 175 rue du Chevaleret, 75013
Paris, France.} \email{picard@math.jussieu.fr}
\author{Karine Tribouley}
\address{Universit\'e Paris Ouest, 92001 Nanterre Cedex, France.} \email{karine.tribouley@u-paris10.fr}

\begin{document}

\begin{abstract}\hspace{0.4cm}
This paper investigates the estimation problem in a regression-type
model. To be able to deal with potential high dimensions, we provide
a procedure called LOL, for Learning Out of Leaders with no
optimization step. LOL is an auto-driven algorithm with two
thresholding steps. A first adaptive thresholding helps to select
leaders among the initial regressors in order to obtain a first
reduction of dimensionality. Then a second thresholding is performed
on the linear regression upon the leaders. The consistency of the
procedure is investigated. Exponential bounds are obtained, leading
to minimax and adaptive results for a wide class of sparse
parameters, with (quasi) no restriction on the number $p$ of
possible regressors.  An extensive computational experiment is
conducted to emphasize the practical good performances of LOL.

\end{abstract}

\vspace{1cm}


%\baselineskip=18 pt
%%%%%%%%%%%%%%%%%%%%%%%%%%%%%%%%%%%%%%%%%%%%%%%%%%%%%%%
\section{Introduction}
%%%%%%%%%%%%%%%%%%%%%%%%%%%%%%%%%%%%%%%%%%%%%%%%%%%%%%%
The general linear model is  considered here, with a particular
focus on cases where the number $p$ of regressors is large compared
to the number $n$ of  observations (although there is no such
restrictions). These kinds of models have today a lot of practical
applications in many areas of science and engineering including
collaborative filtering, machine learning, control, remote sensing,
and computer vision just to name a few of them. Examples in
statistical signal processing and nonparametric estimation include
the recovery of a continuous-time curve or a surface from a finite
number of noisy samples. Other interesting fields of application are
 radiology and biomedical imaging when fewer measurements
 about an image are available compared to the
unknown number of pixels collected. In biostatistics,
 high dimensional problems frequently arise specially in genomic
when gene expression are studied given a huge number of initial
genes compared to a relatively low number of  observations.

A considerable amount of work has been produced in this domain in
the last years, which has been a large source of inspiration for
this paper:  algorithms coming from the learning framework
 \cite{MR2387964}, \cite{MR2249856}, \cite{MR2327596}, \cite{MR2335129}), as well as
the extraordinary explosive domain of $\ell_1$ penalties, among many
others \cite{Lasso1}, \cite{dantzig}, \cite{Lasso3},
 \cite{MR2312149}, \cite{MR2397610}, \cite{fanlv}
 and \cite{candesplan}.  See also \cite{lounici-2008} and
 \cite{alquier-2009}.

 The essential motivation of this work is to provide one of the
 simplest procedures
 that achieves, in the same time, good performances.
LOL algorithm (for Learning Out of Leaders) consists in a two steps
thresholding  procedure. As there is no optimization step, it is
important to address the following question: what are the domains where the procedure is competitive
compared to more sophisticated algorithms, especially to algorithms
performing one or two steps $\ell_1$
 minimization ?
One of our aim here is not only to delimit where LOL is competitive but
also to point out where  the simplicity of LOL induces a slight lack of efficiency
 from both a theoretical point of view as from a practical aspect.

Let us start by introducing the ideas of the emergence of LOL
algorithm. This simple procedure can be viewed as an 'explanation'
or as a 'cartoon' of
 $\ell_1$ minimizations. It is well known that
when the regressors are normalized and orthogonal, $\ell_1$
minimization corresponds to  soft thresholding  which itself is
 close to hard thresholding. Hence, it is quite natural to expect
that  thresholding should perform  well, at least in cases not too
far from these orthonormal conditions. It corresponds, as specified
below, to small coherence conditions. A tricky problem occurs when
the regressors are not orthonormal or when the number of regressors
is large. Then, the minimum least squares estimator has a non unique
solution and the solutions are very unstable. This is the heart and
the main difficulty for the $\ell_1$ minimizers or more generally
for all methods based on sparsity assumptions. In order to be
solved, the problem requires essentially -as will be discussed
extensively in the sequel- two types of conditions: sparsity of the
solution and isometry properties for the matrix of regressors. This
is often the part  where the algorithms computation
 cost shows up. Obviously a simple thresholding would
not fit, but the above mentioned conditions can ensure that it is at
least possible to select some regressors and exclude some others.
LOL algorithm solves the difficult problem of the choice of the
regressors in a quite crude way by adaptively selecting  $N$
regressors which are the most correlated to the target: this defines
the first step thresholding of LOL, determining the $N$ leaders. The
number $N$ is chosen using a fine tuning parameter depending on the
coherence. It has to be emphasized that the choice is auto driven in
the algorithm. In a second thresholding step, LOL regresses on the
leaders, then thresholds the estimated coefficients taking into
account the noise of the model.

The properties of LOL are here investigated specifically for the
prediction problem. More precisely, it is established in this paper
that LOL has a prediction
 error which is going to zero in probability with exponential rates.
 These types of results are often called Bahadur type efficiency.
  Although Bahadur efficiency of test and estimation procedures goes
  back to the sixties (see \cite{Bahadur}), it has seen recently
  a revival  in learning theory, where the rates of convergence
  (preferably exponential) of the procedures are investigated and compared to optimality. It is also related to a
common concept in learning theory: The Probably Approximately
Correct (PAC) learning paradigm introduced in \cite{valiant}.

Of course, because of the straightforwardness  of the method, some
loss of efficiency is expected compared to more elaborate and costly
procedures. But even with a loss,
 the limitations of the procedure can bring an interesting
 information on the $\ell_1$ minimizers themselves. From both  theoretical
  and  practical point of view,
    with small coherence,  LOL procedure appears to be as powerful as the best
    known procedures. The exponential rates of convergence match for instance the lower bounds
    obtained in \cite{raskutti2009}.
    This  result is obtained under minimal conditions on the number $p$ of
   potential regressors.
 Also even with a loss in the rate, a positive aspect
 is that the practitioner is informed of the possible instability of the method
  since the coherence can be computed using the observations before any calculation.
 This is notably not the case for usual conditions such as RIP or
 even more abstract ones which are impossible to verify in practice.
An intensive calculation program is performed to show the advantages
 and  limitations of LOL procedure in several practical aspects. The case where
 the regressors are forming a random design matrix with i.i.d. entries is investigated in Section~\ref{section_simu}.
Different laws of the inputs are studied (Gaussian, Uniform,
Bernoulli or Student laws) inducing a specific coherence for the
design matrix.  Several interesting features are also discussed in
this section. Dependent inputs are simulated and an application with
real data is also discussed. The impact of the sparsity and the
undetermination of the regression on the performances of LOL are
studied. A comparison with two others two-step procedures namely
\cite{fanlv} and \cite{candesplan} is also performed. The most
interesting conclusion being that the practical results are even
better and more comforting than the theoretical ones in the sense
that LOL shows good performances, even when the coherence is pretty
high.

To summarize this presentation  and  answer to the question "In what
type of situations should a practitioner prefer to use LOL rather
than other available methods?", our results and our work prove that
when the number of regressors $p$ is very large, and when the
computational aspects of optimization procedures become difficult as
well as the theoretical results uncertain, LOL should be preferred
by a practitioner after ensuring (and this is done by a simple
calculation) that the coherence is not too high. On the other hand
when the coherence is very high, one should probably be suspicious
enough regarding any method...

The paper is organized as follows. In Section~\ref{section_models},
the general model  and the notations  are presented. In
Section~\ref{section_procedures},  LOL is detailed as other
procedures with a $\ell_1$ optimization step; Comparisons with other
procedures are later discussed in Section~\ref{section_comparison}.
In Section~\ref{section_theorie}, after stating the hypotheses
needed on  the model, theoretical results are established. The
practical performances of LOL  are investigated in
Section~\ref{section_simu} and
 the proofs are detailed in Section~\ref{section_proofs}.

%%%%%%%%%%%%%%%%%%%%%%%%%%%%%%%%%%%%%%%%%%%%%%%%%%%%%%%
\section{Model and coherence}\label{section_models}
%%%%%%%%%%%%%%%%%%%%%%%%%%%%%%%%%%%%%%%%%%%%%%%%%%%%%%%
\subsection{General model}
%%%%%%%%%%%%%%%%%%%%%%%%%%%%%%%%%%%%%%%%%%%%%%%%%%%%%%%
In this paper, we observe a pair $(Y, \Phi)\; \in \R^n\times
\R^{n\times p} $ where $\Phi$ is the design matrix and $Y$ a vector
of response variables. These two quantities are linked by the
standard linear model
\begin{align}
 Y&=\Phi\alpha+u+\eps\label{model}
 \end{align}
 where the parameter
  $\alpha\in \R^p$ is the unknown vector to be estimated and
 \begin{itemize}
\item the vector $\eps= (\eps_1,\ldots,\eps_n)^t $
  is a (non observed) vector of random errors. It is assumed to
   be independent Gaussian variables  $N(0,\sigma^2)$ but
 essentially comparable results can be obtained in the case of zero mean
 subgaussian errors
 (see the remark before Lemma \ref{chi2}).
  \item the vector
 $u= (u_1,\ldots,u_n)^t$ is a non observed vector of (possibly) random errors.
 Its amplitude is assumed to be small. The differences between
the two previously described "errors" lies in the fact that the
$\eps_i$'s are centered but unbounded and independent, while the
$u_i$'s are only bounded. The necessity of
 introducing these two types of errors becomes clear in the
 functional regression example.
   \item $\Phi$ is a $n\times p$ known matrix.  This paper
   focuses on the interesting case where $p\gg n$ but it is
   not  necessary. We  assume that $\Phi$ has normalized columns
 (or normalized them) in the following sense:
 \begin{equation}
\frac 1n\sum_{i=1}^n \Phi_{i\ell}^2=1\label{norm}, \quad \forall \; \ell=1\ldots, p.
\end{equation}
\end{itemize}

\subsection{Examples}
An example of such a model occurs when the matrix $\Phi$ is a
 random matrix  composed of $n$ independent and mainly identically distributed
 random vectors of size $p$. The simulation study given in Section \ref{section_simu}
  details the important role played by the distribution of these random vectors.

A second application is the learning (also called functional
regression) model
\begin{equation}
Y_i= f(X_i)+\eps_i, \; i=1\ldots n\label{functreg}
\end{equation}
where $f$ is the parameter of interest. This model is classically
related to the previous one using a dictionary $\D=\{g_l,\; l\le
p\}$ of size $p$,  $\Phi$ becoming then the matrix with general term
$\Phi_{i\ell}=g_\ell(X_i)$. Assuming that $f$ can be reasonably well
approximated using the elements of the dictionary  means that $f$
can be written as $f= \sum_{g\in \D}\alpha_g g +h$ where $h$ is
hopefully small. It becomes clear here that $u_i=h(X_i)$. This case
has been investigated in more details in \cite{KMPT2009}) using an
earlier and less elaborated version of LOL.

%%%%%%%%%%%%%%%%%%%%%%%%%%%%%%%%%%%%%%%%%%%%%%%%%%%%%%%
\subsection{Coherence}
%%%%%%%%%%%%%%%%%%%%%%%%%%%%%%%%%%%%%%%%%%%%%%%%%%%%%%%
In the sequel, the following notations are used. Let $m$ be an
integer and $q>0$, for any $x\in\R^m$,
$$\|x\|_{l^q(m)}\,:= \,\left(\sum_{k=1}^m |x_k|^q\right)^{1/q}$$
 denotes the $l^q(\R^m)-$norm (or quasi norm) and, for any $x\in \R^n$,
$$\|x\|_n^2\,:=\,\frac 1n \|x\|_{{l}_2(n)}^2$$ denotes the quadratic empirical norm.
 We define the following
$p\times p $ Gram matrix as
\begin{equation*} M:=\frac 1n \Phi^t \Phi.
\label{gram}
\end{equation*}
The quantity
\begin{equation*}
\tau_n=\sup_{\ell\not=m }|M_{\ell m}|=\sup_{\ell\not=
m}|\frac1n\sum_{i=1}^n \Phi_{i\ell}\Phi_{im}|
\end{equation*}
is called the coherence of the matrix $M$. Observe that $\tau_n$ is a  quantity directly
 computable from the data. It
is also a crucial quantity because it induces a bound on the size of
 the invertible matrices built with the columns of $M$.
More precisely, fix $0<\nu<1$ and  let $\cI$ be a subset of indices
of $\{1,\ldots,p\}$ with cardinality $m$. Denote $\Phi_{|\cI}$ the
matrix restricted to the columns of $\Phi$ whose indices belong to
$\cI$. If $2\tau_n \leq\nu$, the associated Gram matrix
$$M(\cI):=\frac 1n
\Phi_{|\cI}^t \Phi_{|\cI}$$
 is almost diagonal as soon as $m$ is smaller than $N:=\lfloor\nu/\tau_n\rfloor$ (where $\lfloor\nu/\tau_n\rfloor$ denotes the integer part of $\nu/\tau_n$)
  in the sense that it satisfies
the following so called Restricted Isometry Property (RIP)
\begin{equation}\label{cond}
\forall x\in \R^m,\;     \       \|x\|_{{ l}_2(m)}^2(1-\nu)\le
x^t\M(\cI)x \leq   \|x\|_{{ l}_2(m)}^2(1+\nu) \;. \end{equation}
This proves in particular that the matrix $M(\cI)$ is invertible.
The proof of \eref{cond}
 is simple and can be found together with a discussion on the relations between RIP Property and conditions on the coherence, for instance in \cite{tanner}.
The RIP Property \eref{cond} can be rewritten as follows. The following lemma is a key ingredient of our proofs.
%%%%%%%%%%%%%%%%%%%%%%%%%%%%%%%%%%%%%%%%%%%%%%%%%%%%%%%
\begin{lemma}\label{rho-eucl}Let $\cI$ be a subset
of $\{1,\ldots,p\}$ satisfying $\#(\cI)\leq N$. For any $x\;
\in\R^{\#(\cI)}$, we get
$$(1-\nu)\,\|x\|_{l^2(\#(\cI))}^2\le \|\sum_{\ell\in\cI} x_\ell\;\Phi_{\bullet\ell}
\;\|_{n}^2\le (1+\nu)\,\|x\|_{l^2(\#(\cI))}^2.$$
\end{lemma}

%%%%%%%%%%%%%%%%%%%%%%%%%%%%%%%%%%%%%%%%%%%%%%%%%%%%%%%
\section{Estimation procedures}\label{section_procedures}
%%%%%%%%%%%%%%%%%%%%%%%%%%%%%%%%%%%%%%%%%%%%%%%%%%%%%%%

In this section,  the
estimation of the unknown parameter $\alpha$ using LOL is described first.
Next, a short review on the procedures directly connected to LOL is
proposed.

Once for all, the constant $\nu$ is fixed. This constant is
obviously  related to the precision of LOL main procedure: the
default value here considered is $\nu=0.5$.

%%%%%%%%%%%%%%%%%%%%%%%%%%%%%%%%%%%%%%%%%%%%%%%%%%%%%%%
\subsection{LOL Procedure}
%%%%%%%%%%%%%%%%%%%%%%%%%%%%%%%%%%%%%%%%%%%%%%%%%%%%%%%
As inputs, LOL algorithm requires 4 pieces of information:
\begin{itemize}
\item The observed variable $Y$ and the regression variables
$\Phi=(\Phi_{\bullet 1},\ldots \Phi_{\bullet p})$.
\item The tuning parameters $\lan(1)$ and $\lan(2)$ giving the level
of the thresholds.
\end{itemize}
Observe that the algorithm is adaptive in the sense that no
information on the sparsity of the sequence $\alpha$ is necessary.
An upper bound for the cardinal of the set of  leaders is first
computed: $N \leftarrow \lfloor
 \nu/\tau_n\rfloor$. Thereafter, LOL performs two major steps:
\begin{itemize}
\item Find the leaders by thresholding the "correlation" between
$Y$ and the $\Phi_{\bullet\ell}$'s at level $\lan(1)$. $\cB$
denotes the set of  indices of the selected leaders. The size of
this set is bounded with $N$ by retaining only (when necessary)
the indices with maximal correlations.
\item Regress $Y$ on the leaders
$\Phi_{|\cB}=(\Phi_{\bullet\ell})_{\ell\in\cB}$ and threshold
the result at  level $\lan(2)$.
\end{itemize}
The following pseudocode gives details of the procedure. Note that
there is no step of optimization and no iteration procedure.

\begin{center}
\label{algo}
\begin{tabular}{|l|}\hline
\\
{\bf LOL$(\Phi,Y,\lan(1),\lan(2))$}\\
\\
{\bf Input:} observed data $Y$, regression variables $\Phi$,
tuning parameters $\lan(1),\lan(2)$\\
{\bf Output:} estimated parameters $\widehat{ \alpha^*}$, and
predicted value
$\widehat{Y}$\\
\\ \hline\end{tabular}

{\footnotesize
\begin{tabular}{|lr|}\hline
{\bf STEP 0}& \{Initialize\}\\
$\nu=0.5$&\\
 $\tau_n \leftarrow n^{-1}\,\max_{\ell\not=
m}|\sum_{i=1}^n\Phi_{i\ell}\Phi_{im}|$&\{Compute the
coherence\}\\
 $N \leftarrow \lfloor
 \frac{\nu}{\tau_n}\rfloor$ &\{Compute the
 upper bound \\&for the cardinal of the leaders set\}\\
& \\
{\bf STEP 1} & \{Find the leaders\}\\
For $\ell=1:p$&\\
\hspace{0.5cm}$\widetilde{\alpha_\ell} \leftarrow
 \frac 1n\,\sum_{i=1}^n\Phi_{i\ell}Y_{i}$&\{Compute the
'correlations' \\&between the observations and the regressors
\}\\
\hspace{0.5cm}$\widetilde{\alpha_\ell}^* \leftarrow
\widetilde{\alpha_\ell}\I\{|\widetilde{\alpha_\ell}|\geq
\lan(1)\}$&\{Threshold
\}\\
End(for)&\\
$\cB \leftarrow \{\ell,\widetilde{\alpha_\ell}^*\not =0\}$
&\{Determine the leaders\}\\
If $\#\cB>N$ &\\
\hspace{0.5cm}indices $\leftarrow$ sort$(|\widetilde{\alpha}|)$&\{Sort the 'correlations' of the candidates\}\\
\hspace{0.5cm}$\cB \leftarrow$ indices$[1:N]$
&\{Take the indices associated to the $N-$th largest\}\\
End(if)&\\
&\\
{\bf STEP 2} & \{Regress on the leaders\}\\
$\widehat{\alpha}_{|\cB} \leftarrow
(^t\Phi_{|\cB}\Phi_{|\cB})^{-1}\Phi_{|\cB}Y $& \{Least square estimators\}\\
 $\widehat{\alpha}^*_{|\cB} \leftarrow
\widehat{\alpha}_{|\cB}\,\I\{|\widehat{\alpha}_{|\cB}|\geq \lan(2)\}
$
&\{Threshold\}\\
$\widehat{\alpha}^*_{|\cB^c} \leftarrow 0 $&
\\
$\widehat{Y}\leftarrow \Phi\,\widehat{\alpha}^*$ & \{Find the
predicted value\}\\&\\\hline
\end{tabular}
}
\end{center}

\subsection{Several inspirations}\label{section_inspiration}
%%%%%%%%%%%%%%%%%%%%%%%%%%%%%%%%%%%%%%%%%%%%%%%%%%%%%%%

Although it is impossible to be exhaustive in such a productive
domain, some of the works directly in relation to our construction
are hereafter mentioned. We apologize in advance for all the works
that are not mentioned  but still remain in connection. For a
comprehensive overview, we refer to \cite{FanLv2009}.

Several authors propose procedures to solve  the selection problem
 or the estimation problem in cases where the vector $\alpha$
 has only a small number of non zero components, and (often) when  the design matrix
 $\Phi$ is composed of i.i.d. random vectors:  see among many others \cite{Lasso1},
 \cite{dantzig}, \cite{Lasso3},
 \cite{MR2312149} and \cite{MR2397610}.

A focus is particulary made here on the 2-steps procedures which are
also commonly used, and apparently for a long time, since in 1959
such a procedure is already discussed (see \cite{satter}). In
\cite{dantzig} and \cite{candesplan}, the leaders are selected with
respectively the Dantzig procedure and the Lasso
 procedure. Then, the estimated coefficients are computed using a linear
regression on the leaders. Using an intensive simulation program,
\cite{fanlv} show that it could be unfavorable to use the procedures
Lasso or Dantzig {\it before} the reduction of the dimension. They
also provide a search among leaders called Sure Independence
Screening (SIS) procedure. This procedure is similar to the one
discussed in this paper: The leaders are the $N=\lfloor \gamma_n
n\rfloor$ columns of $\Phi$ with largest correlations to the target
variable $Y$ ($\gamma_n$ is a tuning sequence tending to zero). This
step is followed with a subsequent estimation procedure using
Dantzig or Lasso. All these methods show a  higher complexity
compared to LOL.

LOL procedure can also be connected to the family of Orthogonal
Matching Pursuit algorithms as well as in general to  the Greedy
Algorithms. For this interesting literature, we refer among others
to \cite{needell-2008}, \cite{Tropp-Guilbert}, \cite{MR2387964}. The
main advantage of LOL compared to this kind of algorithms is that
there is no iterative search of the leaders. All the leaders are
selected in one shot and the procedure stops just after the second
step. Moreover, convergence results are almost as good as those
procedures in many situations.

%%%%%%%%%%%%%%%%%%%%%%%%%%%%%%%%%%%%%%%%%%%%%%%%%%%%%%%
\section{Main theoretical results}\label{section_theorie}
%%%%%%%%%%%%%%%%%%%%%%%%%%%%%%%%%%%%%%%%%%%%%%%%%%%%%%%

This section states the theoretical results of LOL procedure. The
measures of performances used in the theorems are first presented,
 then the assumptions on the set of parameters $\alpha$ are given.

\subsection{Loss fonction} \label{perfdefinition}
%%%%%%%%%%%%%%%%%%%%%%%%%%%%%%%%%%%%%%%%%%%%%%%%%%%%%%%

Let us define the following loss function to measure the difference
between the true value $\alpha \; \in \R^p$ and the result $\hat
\alpha^*$ computed by LOL. Denote $\Phi_{i\bullet}$ the $i-$th line
of the matrix $\Phi$ and recall that the $i-$th observation is given
by the model:
$$Y_i=\Phi_{i\bullet}\alpha+u_i+\epsilon_i.
$$
The predicted $i-$th observation is
$\widehat{Y_i}=\Phi_{i\bullet}\hat\alpha^*$. The criterium of
performance is defined by the empirical quadratic distance between
the predicted variables  and their expected values.
\begin{eqnarray*}\label{distance}d(\hat\alpha^*,\alpha)^2=\frac1n\sum_{i=1}^n
\left(\widehat{Y_i}-\E Y_i\right)^2
=\frac1n\sum_{i=1}^n\left(\sum_{\ell=1}^p(\hat\alpha_\ell^*-\alpha_\ell)\Phi_{i\ell}+u_i\right)^2
\end{eqnarray*}
which can be rewritten using the empirical norm
$$d(\hat\alpha^*,\alpha):=\|\sum_{\ell=1}^p(\hat\alpha_\ell^*-\alpha_\ell)\Phi_{\bullet\ell}
+u_\bullet\|_{n}.$$ Observe that  the considered loss is the usual
error of prediction
$$d(\hat\alpha^*,\alpha)= \|\Phi(\hat\alpha^*-\alpha)\|_{n}$$
when the 'errors' $u_i$'s are all zero.

%%%%%%%%%%%%%%%%%%%%%%%%%%%%%%%%%%%%%%%%%%%%%%%
\subsection{Bahadur-type efficiency}

Our measure of performance  is issued from
the Bahadur efficiency of test and estimation procedures and is
defined for any tolerance $\eta>0$ as
 \beqn {AC}_n(LOL,\eta,\alpha)=P\left(
d(\hat\alpha^*,\alpha)>\eta\right). \eeqn
This quantity  measures a  confidence
   that the estimator  is accurate to the tolerance $\eta$ if the true point is $\alpha$.
   We also define and consider uniform confidence over a class $\Theta$ :
 \beqn {AC}_n(LOL,\eta,\Theta)=\sup_{\alpha\in\Theta}P\left(
d(\hat\alpha^*,\alpha)>\eta\right). \eeqn
This quantity has been studied for instance in  \cite{dkpt} in the learning framework.
 In most examples, it is proved that there exist a phase transition and a critical value
  $\eta_n$ depending on $n$ and $\Theta$ such that
  ${ AC}_n(\hf,\eta,\Theta)$ decreases exponentially for any $\eta>\eta_n$.
More precisely, in terms of lower bound -but similar bounds are also
valid in terms of upper bounds-, it is proved in \cite{dkpt} that
\beqn \label{acc2} \inf_{\hf}{ AC}_n(\hf,\eta,\Theta)\ge C\sqrt
{\bar N(\Theta,\eta)}e^{-cn\eta^2}, \eeqn
where $\bar N(\Theta,\eta)$ is the tight entropy  analogue of the
Sobolev covering numbers. On this expression,  $\eta_n$ appears
quite convincingly as  a turning point after which the exponential
term dominates the entropy term. Observe that the critical value
$\eta_n$ is essential since it
   yields bounds for $\sup_{\alpha\in \Theta}E_\alpha d(\hat\alpha^*,\alpha)$ which is
  another (more standard) measure of performance of the procedure.
   The results in \cite{dkpt} are obtained in the learning framework; however
   identical bounds can easily be expected in the setting
   \eref{model} of this paper, as results obtained in \cite{raskutti2009}.
    This is discussed in more details in the sequel since similar bounds are obtained for LOL
     with sparsity constraints defined below.

%%%%%%%%%%%%%%%%%%%%%%%%%%%%%%%%%%%%%%%%%%%%%%%%%%%%%%
\subsection{Performances of the procedure LOL.  $l_q$ ball constraints}
%%%%%%%%%%%%%%%%%%%%%%%%%%%%%%%%%%%%%%%%%%%%%%%%%%%%
 In this part, we consider the following sparsity constraint
 \begin{align*}
\mbox{for }q\in (0,1],\quad B_q(M):=\{\alpha\in \R^p,\;
\|\alpha\|_{l^q(p)}\le M\}
\end{align*}
or
 \begin{align*}
\mbox{for }q=0,\quad B_0(S,M):=\{\alpha\in \R^p,\;
\sum_{j=1}^pI\{|\alpha_j\not=0\}\le S, \;
  \|\alpha\|_{l^1(p)}\le M\}.
\end{align*}
 \begin{theorem}\label{lqresults} Let $M>0$ and fix $\nu$ in $]0,1[$.
 Assume that there exists a positive constant $c^\prime$ such that $p\leq \exp(c^\prime n
 )$.
  Suppose there exist positive constants $c,c_0$ such that
 \beqn
 \tau_n\le c\sqrt{\frac{\log p}n} \quad\mbox{ and }\quad  \sup_{i=1,\ldots,n}|u_i|\le c_0\;\sqrt{\frac {1}n}.
 \eeqn
Let us  choose the
thresholds $\lan(1)$ and $\lan(2)$ such that
$$
\lan(2)=  T_3\sqrt{\frac{\log{p}}{n}}\quad\mbox{ and }\quad \lan(1)=
T_4\sqrt{\frac{\log{p}}{n}}
$$ for $ T_4\geq T_1\vee T_2\,c\vee T_3>0$ where
$$
T_1=\left(64\sigma\vee1\vee\frac{2}{(1-\nu)\sigma}\right)\mbox{ and
} T_2= \left(6M\vee \frac{(4M+3c_0)}{12\sigma}\right).
$$
Then,
 there exist positive constants $D$ and $\gamma$ depending on $\nu,c,c^\prime,c_0,T_3,T_4$ such that
\begin{equation*}
\sup_{ \alpha\in\; B_q(M)}\P\left(d(\hat\alpha^*,\alpha)>\eta\right)\le
 \left \{
\begin{array}{lll}
4e^{-\gamma n\eta^2} & \mbox{ for }&\eta^2\ge D\;\left(\frac{\log p}n\right)^{1-q/2}\\
&&\\
 1& \mbox{ for }& \eta^2\le D\;\left(\frac{\log p}n\right)^{1-q/2}
\end{array}\right.
\end{equation*}
and
\begin{equation*}
\sup_{  \alpha\in\;  B_0(S,M)}\P\left(d(\hat\alpha^*,\alpha)>\eta\right)\le
 \left \{
\begin{array}{lll}
4e^{-\gamma n\eta^2} & \mbox{ for }&\eta^2\ge D\;\frac{S\,{\log p}}n\\
&&\\
 1& \mbox{ for }& \eta^2\le D\;\frac{S\,\log{ p}}n
\end{array}\right.
\end{equation*}
for any $S<  \nu/\tau_n$.
\end{theorem}
We immediately deduce the following bound for the usual expected error
\begin{corollary}\label{corollary_minmax} For $r\geq 1$ arbitrary, under the same
  assumptions as in Theorem \ref{lqresults}, we have
$$ \sup_{ B_q(M)}\E d(\hat\alpha^*,\alpha)^r\le D^\prime
\left(\frac {\log p}n\right)^{(1-q/2)r/2}$$ for some positive
constant $D^\prime$ depending on $\nu,c,c^\prime,c_0,T_3,T_4$, as
well as
$$ \sup_{ B_0(S,M)}\E d(\hat\alpha^*,\alpha)^r\le D^\prime
\left(\frac {S\log p}n\right)^{r/2}$$ for any $S< \nu/\tau_n$.
\end{corollary}

%%%%%%%%%%%%%%%%%%%%%%%%%%%%%%%%%%%%%%%%%%%%%%%%%%%%%%%
\subsection{Performances of LOL procedure. MaxiSet point of view}
%%%%%%%%%%%%%%%%%%%%%%%%%%%%%%%%%%%%%%%%%%%%%%%%%%%%%%%

In this section we develop a slightly different point of view issued
from the maxiset theory (see for instance \cite{kp-b}). More
precisely our aim is to evaluate the quality of our algorithm when
the coherence and thresholding tuning constants are given (fixed).
Especially, it means that we do not assume in this part that the
coherence satisfies $\tau_n\leq O(\sqrt{\log{p}/n})$. For that, we
consider a set $V(S,M)$ of parameters $\alpha$ depending on these
constants $M,S>0$ and prove that the right exponential decreasing of
the confidence is achieved on this set. The phase transition
$\eta_n$ is depending on the tuning constants and of the coherence
in the following way $$\eta_n^2 = O(\frac { S\log p}n\vee
S\tau_n^2).$$

Observe that we do not prove that the set $V(S,M)$ is exactly the
maxiset of the method (considered in terms of $\eta_n$) since we do
not not prove that it is the largest set with the phase transition
$\eta_n$. However the following theorem  reflects quite extensively
the theoretical behavior of LOL, even in case of deterioration due
to a high coherence or  a bad choice of the thresholds. Notice also
that Theorem \ref{lqresults} is a quite easy consequence of Theorem
\ref{merceradapt}.

Let us now define the set $V(S,M)$ by the following sparsity constraints.
 There exist  $S\le \lfloor\nu/\tau_n\rfloor$  and constants
${M},\; c_0,\; c_1,\,c_2$, such that the vector $\alpha\in \R^p$
satisfies the following conditions
\begin{align}
 \|\alpha\|_{l^1(p)}\le { M}, \label{l1}
 \end{align}
\begin{align}
 \#\left\{\ell\in\{1,\ldots,p\},\;|\alpha_\ell|\ge \lan(2)/2\right\}\;\le
 {S} \label{weak}
 \end{align}
\begin{align} \sum_{(\ell)>N}|\alpha_{(\ell)}|\le c_1 \;
\left(\frac {S\log p}{n\tau_n}\right)^{1/2}\label{l1q}
 \end{align}
\begin{align}\sum_{\ell=1}^p|\alpha_\ell|^2\;I\{|\alpha_\ell |\le 2\lambda_n(1)\}\le
c_2^2\;\frac {S\log p}n\label{l2q}
 \end{align}
Recall that $(\alpha_{(\ell\,)})$ is the ordered sequence (for the
modulus) $|\alpha_{(1)}|\geq
|\alpha_{(2)}|\geq\ldots|\alpha_{(p)}|$. For $S,M>0$,  $V(S,M)$
denotes the class of models of type \eref{model} satisfying the
sparsity conditions  \eref{l1}, \eref{weak}, \eref{l1q}, \eref{l2q}.
Note that we emphasize in the notation of the set $V(S,M)$ the
constants $S$ and $M$, while the set is depending on other
additional constants, since these two constants play a crucial role.

  \begin{theorem}\label{merceradapt} Let $S,M>0$ and fix $\nu$ in $]0,1[$.
The
thresholds $\lan(1)$ and $\lan(2)$ are chosen such that
$$
 \lan(1)\geq
\left(T_{1}\left(\frac{\log{p}}{n}\right)^{1/2}\vee T_{2} \,\tau_n
\right)\quad \mbox{ and }\quad \lan(2)\leq \lan(1)
$$
$$
T_1=\left(64\sigma\vee1\vee\frac{2}{(1-\nu)\sigma}\right)\mbox{ and
} T_2= \left(6M\vee \frac{(4M+3c_0)}{12\sigma}\right).
$$
 Then, if in addition we have
 \begin{align}
 \sup_{i=1,\ldots,n}|u_i|\le c_0\;\left(\frac {S}n\right)^{1/2} \label{h}
 \end{align}
 there exist positive constants $D$ and $\gamma$ depending on $\nu,\sigma^2,M,c_0,c_1,c_2$, such that
\begin{equation}
\sup_{\alpha\in\; V(S,M)}\P\left(d(\hat\alpha^*,\alpha)>\eta\right)\le
 \left \{
\begin{array}{lll}
4e^{-\gamma n\eta^2} & \mbox{ for }&\eta^2\ge D\;\left(\frac { S\log p}n\vee\ S\tau_n^2\right),\\
&&\\
 1& \mbox{ for }& \eta^2\le D\;\left(\frac { S\log p}n\vee S\tau_n^2\right)
\end{array}\right.
\end{equation}
\end{theorem}

 For a sake of completeness,
the constants $D$ and $\gamma$ are precisely given at the end of the
proof of Theorem \ref{merceradapt}.  However, it is obvious that the
constants provided here are not optimal: for instance in the proof,
in order to avoid unnecessary technicalities, most of the events are
divided as if they had an equal importance, leading to constants
which are each time divided by 2. Obviously there is some place for
improvement at any of these stages.

 An elementary consequence of Theorem \ref{merceradapt} is
  the following corollary which details
  the behavior of the expectation of $d(\hat\alpha^*,\alpha)$.
   Notice also that we did not give here explicit oracle inequalities,
   which however could be derived from the proof of Theorem \ref{merceradapt}.

 \begin{corollary}\label{corollary_maxi}For $r\geq 1$ arbitrary, under the same
  assumptions as in Theorem \ref{merceradapt}, we get
$$ \sup_{ V(S,M)}\E d(\hat\alpha^*,\alpha)^r\le D^\prime
\left(\frac { S\log{p}}n\vee S\tau_n^2\right)^{r/2}$$ for some
positive constant $D^\prime$ depending on
$\nu,\sigma^2,M,c_0,c_1,c_2$ and $r$.
\end{corollary}

%%%%%%%%%%%%%%%%%%%%%%%%%%%%%%%%%%%%%%%%%%%%%%%%%%%%%%%
\section{Remarks and Comparisons}
\label{section_comparison}
%%%%%%%%%%%%%%%%%%%%%%%%%%%%%%%%%%%%%%%%%%%%%%%%%%%%%%%
\subsection{Results under $l_q$ constraints} It is important to
discuss the relations of the results in Theorem \ref{lqresults} with
\cite{raskutti2009} which provides minimax bounds in a setting close
to ours. Their results basically concern exponential inequalities
(as ours) but they are only interested in the case $\eta=\eta_n$ for
which they prove upper and lower bounds. If we compare our results
to theirs, we find that LOL is exactly minimax for any $q$ in
$(0,1]$, with even a better precision since we prove the exponential
inequality for any $\eta$. In the case $q=0$, we have a slight
logarithmic loss. Notice that we also need a bound on
$\|\alpha\|_{l^1(p)}$. We do not know if this is due to our proof or
specific to the method.

\subsection{Ultra high dimensions}
One main advantage of LOL  is that it is  really designed for
 very large dimensions. As seen in the results
 of Theorem \ref{lqresults} and  Theorem \ref{merceradapt},
 no limitation on $p$ is required except $p\leq \exp( cn)$ (in fact this is only needed in Theorem \ref{lqresults}).
  Notice that, if this condition is not satisfied, not any algorithm is  convergent
  as proved by the lower bound of \cite{raskutti2009}.
  Moreover, the fact that the algorithm has no optimization step is a
  serious advantage when $p$ becomes large.

\subsection{Adaptation}
Our theoretical results are provided under conditions on  the tuning
quantities $\lan(1)$ and $\lan(2)$. The  default values issued from
the theoretical results are the
 following
\begin{align*} \lan^*(1)=\lan^*(2)= \left(T_{1}\sqrt{\frac{\log{p}}{n}}\vee T_{2}
\,\tau_n \right).
\end{align*}
It is a consequence of  Theorem \ref{lqresults} that  LOL associated
with $\lan^*(1)$ and $\lan^*(2)$ is adaptive over all the sets $
B_q(M)$ and $B_0(S,M)$, with respect to the parameters $q$ and $S$.
These default values behave also reasonably well in practice.
However, they require a fine tuning of the constants $T_{1}$ and
$T_{2}$ which is proposed in a slightly more subtle way in the
simulation part (see Section \ref{section_simu}).

\subsection{Coherence condition}
As can be seen in Theorem \ref{lqresults}, LOL  is minimax under a
condition on the coherence of the type $\tau_n\le c\sqrt{\log{
p}/n}$.
 This condition is
verified with overwhelming probability for instance when the entries
of the matrix $\Phi$ are independent and identically random
variables with a sub-gaussian common distribution. In Section
\ref{section_simu}, we precisely investigate the behavior of  LOL
 when this hypothesis is disturbed. This bound is
 generally stronger as a condition compared to
 other ones given in the literature such as the RIP condition, or weaker ones.
 However, as explained in the sequel, these other conditions
 are often impossible to verify on the data. We consider as
 a benefit that the procedure is giving with $\tau_n$ an indication
 of a potential misbehavior. Besides, Theorem \ref{merceradapt} details
  the behavior of the algorithm when this condition is not verified.

\subsection{Comparison with some existing algorithms}
 As mentioned in the previous section, LOL finds its
inspiration in the learning framework,  especially in
\cite{MR2387964},\- \cite{MR2249856},\-
\cite{MR2327596},\-\cite{MR2335129}. In all these papers,
consistency results are obtained under fewer assumptions but with no
exponential bounds and a higher cost in implementation. Again in the
learning context, \cite{MR2336417} provides optimal critical value
$\eta_n$ as well as exponential bounds with  fewer assumptions since
there is no coherence restriction. However, the procedure is very
difficult to implement for large values of $p$ and $n$ ($N$-$P$
hard).

In \cite{fanlv}, it is assumed that there exists $\kappa>0$ such
that
$$\min_{\ell\in{\mathcal
I}^*}|\alpha_{\ell}|\geq O(n^{-\kappa})$$ where ${\mathcal
I}^*=\{\ell,\alpha_\ell\not = 0\}$. The model under consideration is
ultra high dimensioned: $p\leq\exp( cn^{\xi})$ for $c,\xi>0$ with
the restriction $\xi<1-2\kappa<1$.  The procedure SIS-D (SIS
followed by Dantzig)
 is shown to be asymptotically consistent in the sense that, with large
 probability, we have
$$
\sum_{\ell=1}^p(\hat\alpha_\ell^{SIS-D}-\alpha_\ell)^2\leq C\,\eta_n
 $$
where $C$ is a constant depending on the restricted orthogonality
constant, but the order of  the convergence is not given. A
practical drawback is that the tuning sequence $\gamma_n$ is not
auto driven since it has to verify $\gamma_n=O(n^{-\theta})$ for
$\theta<1-2\kappa-\tau$ for some $\tau$ linked to the largest
eigenvalue of the covariance matrix of the regressors. Notice that
another tuning parameter $\lambda_n$ has also to be chosen in the
Dantzig step.

 In \cite{bunea} and \cite{MR2397610}, the size $p$ grows polynomially
 with the sample size $n$. It is assumed that
\begin{equation*}
\sup_{\ell\in {\mathcal I}^*,\;m\not\in {\mathcal I}^*
}\frac1n\sum_{i=1}^n |\Phi_{i\ell}\Phi_{im}|\leq O(S^{-1}).
\end{equation*}
 which appears to
be a weaker condition on the coherence than ours. However this
condition is impossible to verify on the data since ${\mathcal
I}^*,S$ are unknown. An exponential bound is established for $
P\left(\sum_{\ell=1}^p|\hat
\alpha_\ell-\alpha_\ell|>\sqrt{S}\;\tilde\eta\right)$ when
$\tilde\eta\geq \sqrt{S\log{p}/n}$ corresponding to ours critical
value $\eta_n$. This result is comparable to ours but focuses on the
error due to the estimation of the parameter $\alpha$ instead of the
prediction error.
 In \cite{candesplan}, the condition on the coherence  $\tau_n\leq O((\log{p})^{-1})$
 is  generally lighter except for very large $p$ but no  exponential bounds are provided: it is  proved that
$P\left(d(\hat\alpha_\ell,\alpha_\ell)> \eta\right) $ is tending
to zero as $O\left(p^{-2\log{2}}\right)$ for $ \eta\geq
\sqrt{S\log{p}/n}$.

%%%%%%%%%%%%%%%%%%%%%%%%%%%%%%%%%%%%%%%%%%%%%%%%%%%%%%%
\section{Practical results}\label{section_simu}
%%%%%%%%%%%%%%%%%%%%%%%%%%%%%%%%%%%%%%%%%%%%%%%%%%%%%%%

In this section, an extensive computational study is conducted using
LOL. The performances of LOL are
 studied over various ranges of level of indeterminacy
$\delta=1-n/p$ and of sparsity rates $\rho=S/n$ (see
\cite{malekidonoho2009}). The influence of the design matrix is
investigated: more precisely, as we  consider random matrices, we
study the role of the distribution  for the design matrix $\Phi$ as
well as the nature of dependency between the inputs. This study is
performed on simulations and an application with real data is
presented. Our procedure is finally compared to some others well
known two-step procedures.

%%%%%%%%%%%%%%%%%%%%%%%%%%%%%%%%%%%%%%%%%%%%%%%%%%%%%%%
\subsection{Experimental design}

 The design matrix $\Phi$ considered in this study is generally of
  random type (except in the example of real data) and
   is mostly built  on $n\times p$ independent and identically distributed
inputs. Different distributions such as Gaussian,
Uniform, Bernoulli, or Student laws are considered. We also investigate the influence of the dependency on the procedure.   It is important to stress that all
the above mentioned parameters $p$, $n$, dependency and different type of laws
yield different values of the coherence $\tau_n$ and consequently
different behaviors of the procedure. Each column vector of $\Phi$
is normalized to have unit norms. Given $\Phi$, the target
observations are  $Y=\Phi \alpha + \epsilon$ for $\epsilon$  i.i.d.
variables with a normal distribution $N(0,\sigma^2_{\epsilon})$,
$\sigma_{\epsilon}$ chosen such that the
 signal over noise ratio (SNR) is in most studies SNR=5. When specifies,
 SNR varies  from $SNR=10$ to $SNR=2$. The vector of parameters
 $\alpha$ is simulated as follows: all coordinates are zero except
 $S$ non zero coordinates with
$\alpha_{\ell}=(-1)^ b|z|,\; \ell=1,...,S $ where $b$ is
drawn from a Bernoulli distribution with parameter $0.5$ and $z$
from a $N(2,1)$ (see \cite{fanlv}).

To evaluate the quality of LOL,  the relative $l_2$  error of
prediction $E_Y=\|Y - \hat Y\|_2^2/ \|Y\|_2^2$ is computed on the
target $Y$. The sparsity $S$ is estimated by the cardinal of $
\cL=\{\ell=1,\ldots,p,\widehat\alpha^*\not=0\}$ where
$\widehat\alpha^*$ is provided by LOL. All these quantities are
computed by averaging each estimation result  over $K$ replications
of the experiment ($K=200$).

%%%%%%%%%%%%%%%%%%%%%%%%%%%%%%%%%%%%%%%%%%%%%%%%%%%%%%%
\subsection{Algorithm}
The parameters $\lambda_n(1)$ and $\lambda_n(2)$
 are critical values  quite hard to tune practically
 because they depend on constants which are not optimized and may be unavailable in practice (such as the constant $M$ -see the theoretical results-).
Let us explain how we proceed in this study to adaptively determine the
thresholds.

Since the first threshold $\lambda_n(1)$ is used to select the
leaders, our aim is to split the set of
"correlations" $\{K_\ell,\;\ell=1,\ldots,p\}$,
 into  two clusters in such a way that the leaders are forming one of the two clusters.
The sparsity assumption suggests that the law of
the correlations (in absolute value) should be a mixture of two distributions:  one for the leaders (high correlations- positive mean) and one
 for the others (very small correlations- zero mean).
 The frontier between the clusters is then chosen by minimizing the variance between classes after adjusting the  absolute
 value of the correlations into the  two classes  described above
 ( see also \cite{KMPT2009}).

The same procedure is used to threshold adaptively the estimated
coefficients $\widehat\alpha_\ell$ obtained by linear regression on
the leaders. Again the distribution of the
$\widehat\alpha_\ell$ provides two clusters: one cluster associated
to the largest coefficients (in absolute value) corresponding to the
non zero coefficients  and one cluster composed of coefficients
close to zero, which should not be involved in the model. The
frontier between the two clusters, which defines $\lambda_n(2)$, is
again computed by minimizing the deviance between the two classes of
regression coefficients.

Finally, an additional improvement for  LOL is provided. It
generally more efficient to perform a second regression using  the
final set of selected predictors involved in the model: the
estimators of the (non zero) coefficients are then slightly more
accurate. This updating procedure is denoted LOL$^+$ in the sequel.

%%%%%%%%%%%%%%%%%%%%%%%%%%%%%%%%%%%%%%%%%%%%%%%%%%%%%%%
\subsection{Results with i.i.d. gaussian design matrices}
\begin{figure}
\makebox{\includegraphics[width=6cm,height=6cm]{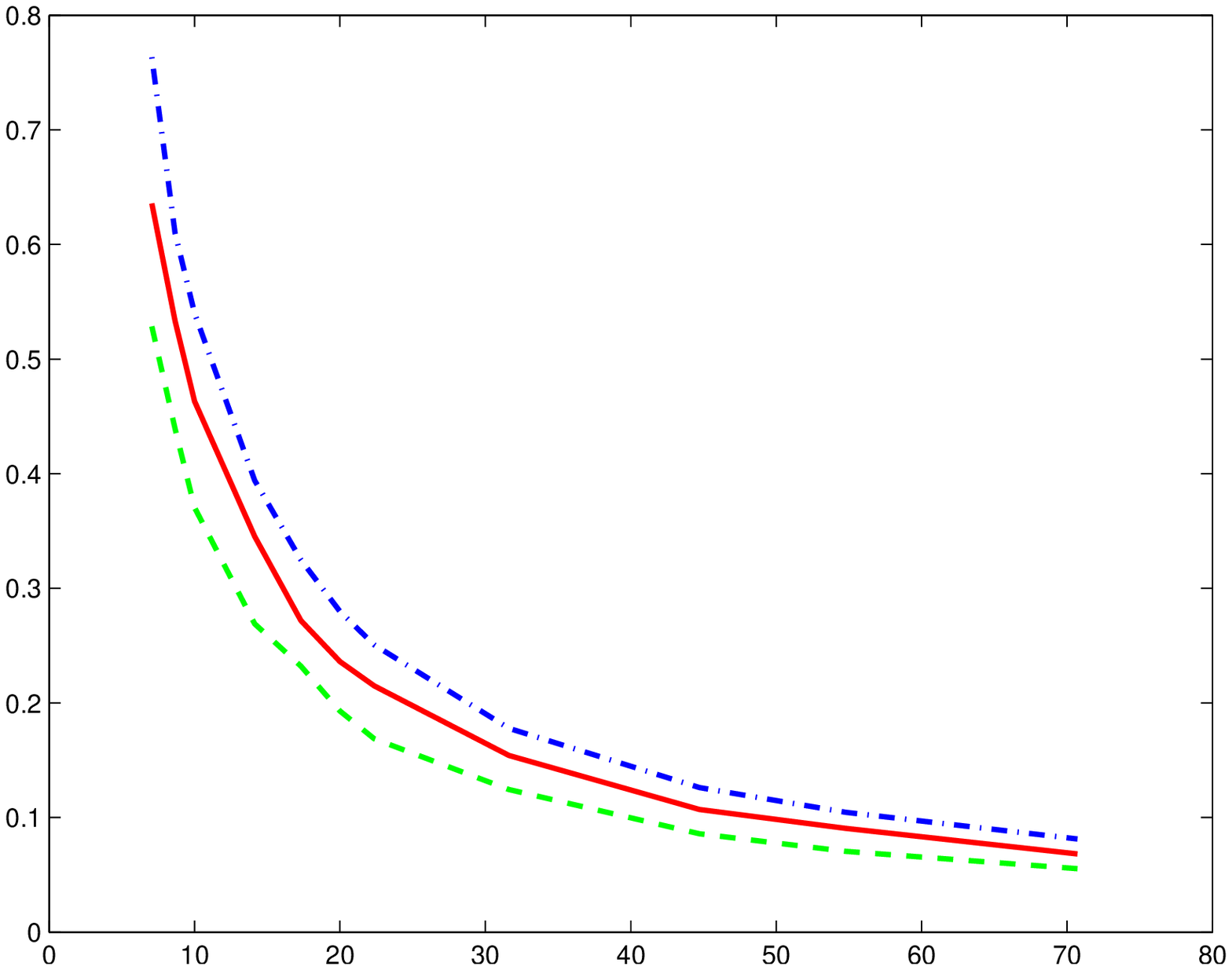}
\includegraphics[width=6cm,height=6cm]{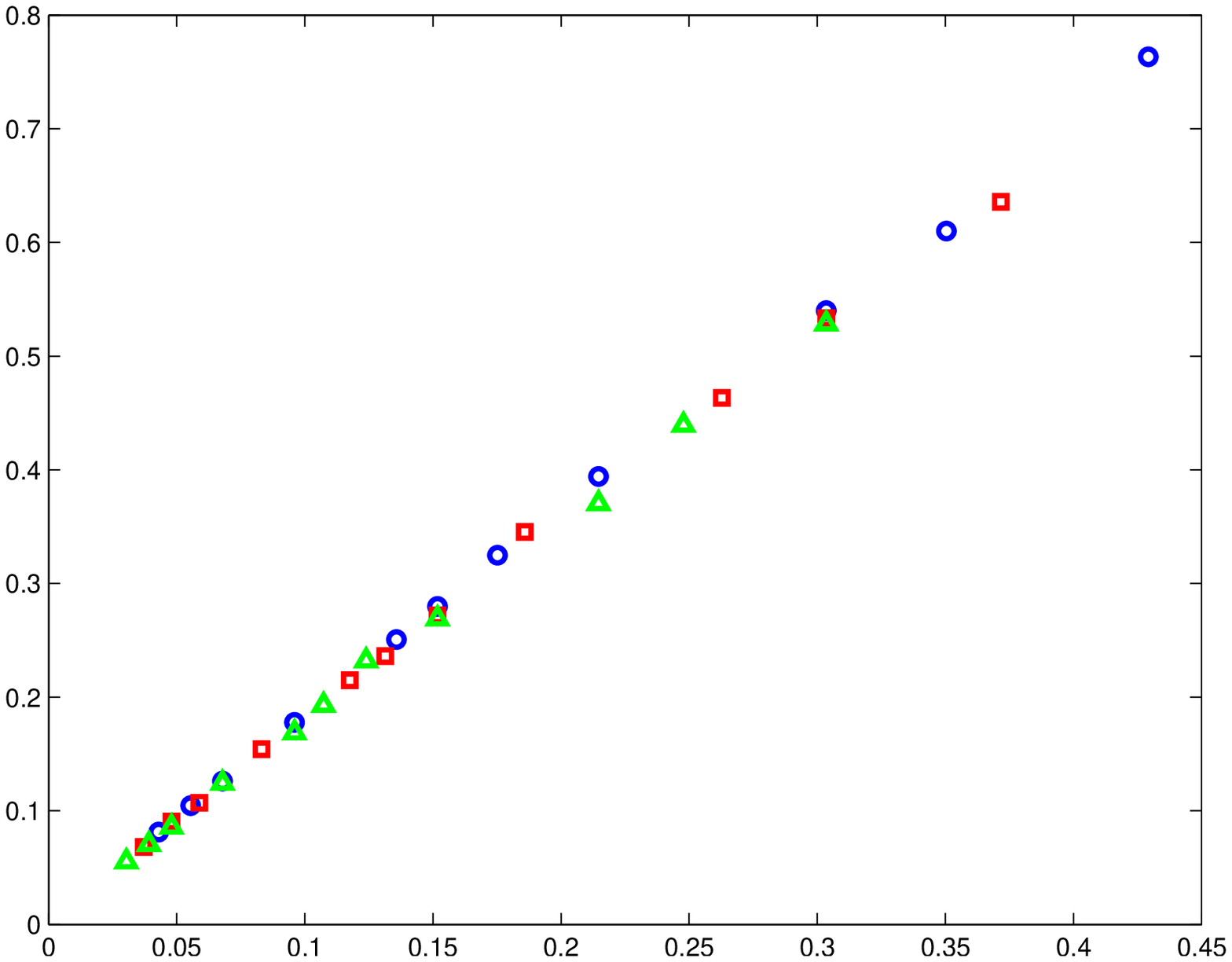}}
\caption{$Y-$axis:  Coherence $\tau_n$.  $X-$axis: $\sqrt{n}$ (left)
 or $\sqrt{\frac{log(p)}{n}}$ (right) for $p=100$ (dashdot line or
triangle -green), $p=1000$ (solid line or square -red), $p=10000$
(dash line or circle -blue). $K=500$ \label{cohstudy}}
\end{figure}

The design matrix $\Phi$ is first defined with i.i.d. gaussian
variables. Figure $\ref{cohstudy}$ (left) shows the evolution of the
empirical coherence $\tau_n$  function of $\sqrt{n}$ for $p=100$,
$1000$, or $10000$. Each coherence  shown in the graph is the
average of $K=500$ coherence values computed for different  $\Phi$
matrix simulated at random over the $K$  replications. As the number
of observations increases to $n=5000$ ($\sqrt{5000} \simeq 70.7$),
the coherence tends to be quite small ($\tau_n=0.1$) independently
of the number of variables $p$. For a small number of observations,
the coherence takes pretty high values, much higher as the number of
predictors increases. For example, for $n=250$ ($\sqrt{250}\simeq
15.8$)
$$ p=100\mapsto \tau_n=0.25,\quad p=100 \mapsto \tau_n=0.30,
\quad p=1000 \mapsto \tau_n=0.35.$$  A difference of $15\%$ is
observed between the coherences computed for $p=1000$ and $p=100$,
or $p=1000$ and $p=10000$. Figure $\ref{cohstudy}$ (right) shows the
evolution of the coherence as a function of $\sqrt{log(p)/n}$ which
allows to compute the constant $c$ introduced in Theorem
$\ref{lqresults}$.

Since we are interested by quantifying the performances of LOL
 in an overwhelming majority of cases ($n$, $p$ varying),
the impact of the level of indeterminacy and of the sparsity rate
are studied: $\delta$ is varying from $0$ to $0.9$ by $0.05$ step
and $\rho$ is varying from $0.01$ to $0.16$ by 20 steps. We fixe
$p=1000$ and $n=250$ for this specific study.

\begin{figure}
\makebox{\includegraphics[width=10cm,height=10cm]{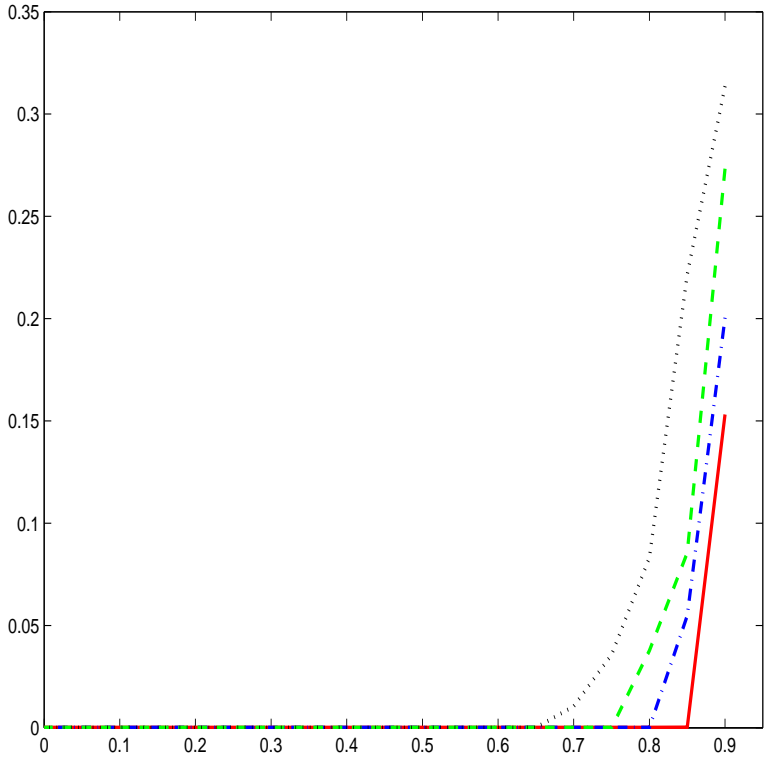}}
\caption{$X-$axis: indeterminacy level $\delta$, $Y-$axis: relative
prediction error. $S=10$ (solid line-red); $S=12$ (dashdot
line-blue); $S=15$ (dash line-green); $S=20$ (dot line-black).
$SNR=5$. \label{LolEY2delta}}
\end{figure}

{\bf Influence of the indeterminacy level:} Figure \ref{LolEY2delta}
studies the performances of LOL when the indeterminacy level is
varying ($p=1000$ fixed, $n$ varying), for different sparsity values
($S=10, 12, 15, 20$). The error of prediction $E_Y$ increases
continuously with the indeterminacy $\delta$, as the number of
 observations decreases compared to the number of
variables. For a given value of $\delta$, $E_Y$ decreases as the
sparsity does.  For $\delta \leq 0.75$, the prediction error is
weak, below $5\%$. When the number of available observations is at
least higher than half of the number of potential predictors
($\delta<0.5$), the prediction  error is negligible:  the quality of
LOL is in this case exceptionally good. For a given number of
observations and potential predictors, the prediction is more
accurate as the sparsity rate decreases. For a fixed number of
observations, regarding the joint values of both indeterminacy and
sparsity parameters, the errors tends to be null as $\delta$ and/or
$\rho$ are decreasing.

\begin{figure}
\makebox{\includegraphics[width=10cm,height=10cm]{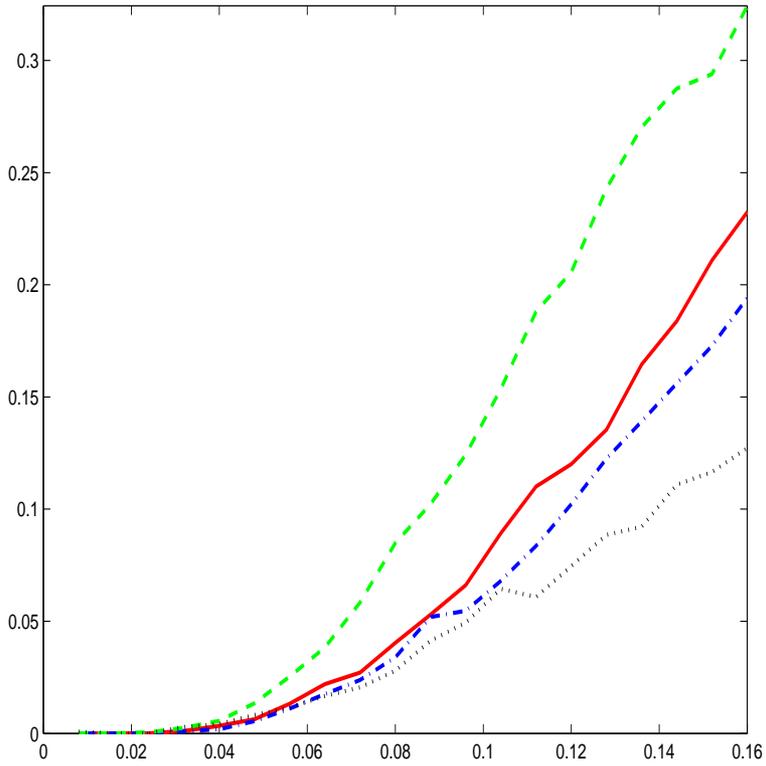}}
\caption{ $X-$axis: sparsity rate $\rho$, $Y-$axis: relative
prediction error. $\delta=0.4$ (dot line-black); $\delta=0.7$
(dashdot line-blue); $\delta=0.75$ (solid line-red); $\delta=0.875$,
(dashed line-green). $SNR=5$. \label{LolEY2rho}}
\end{figure}

{\bf Influence of the sparsity rate:} Figure \ref{LolEY2rho}
illustrates the performances of LOL for prediction when the sparsity
rate is varying for four levels of indeterminacy ($\delta=0.4, 0.7,
0.75, 0.875$). For small values of the sparsity rate ($\rho \leq
5\%$), the prediction error is very good (less than $5\%$). For an
extreme level of sparsity ($\rho \leq 2\%$), the performances are
excellent. As observed before, for a given sparsity rate value, the
performances are improved as the indeterminacy level is decreasing.

\begin{figure}
\makebox{\includegraphics[width=10cm,height=10cm]{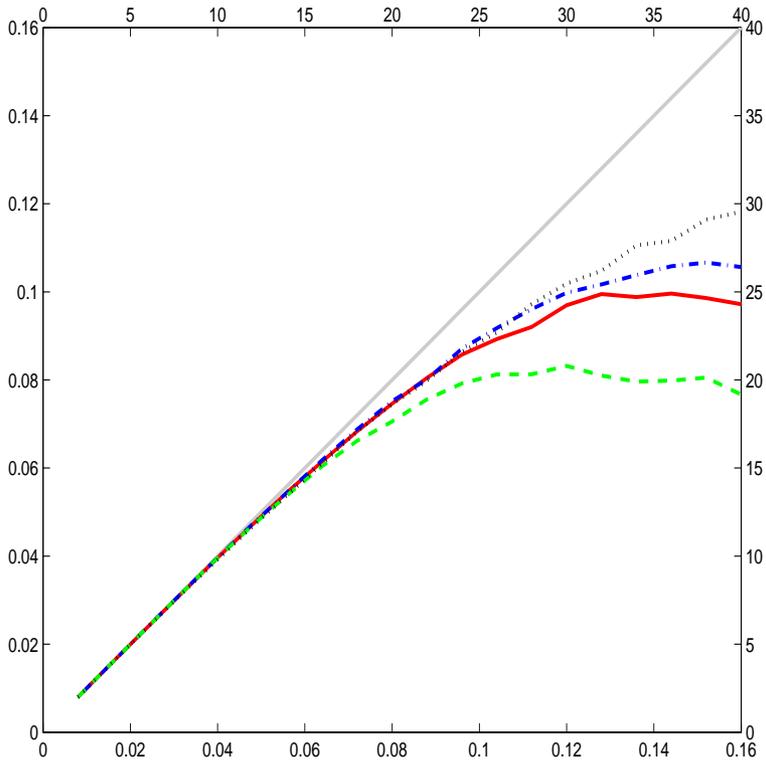}}
\caption{ LOL Sparsity Estimation ( $\rho$: bottom, left; $S$:
right, top). $\delta=0.875$ (dashed line-green); $\delta=0.75$
(solid line-red); $\delta=0.7$ (dashdot line-blue); $\delta=0.4$
(dot line-black). \label{LolRhoEst}}
\end{figure}

{\bf Estimation of the Sparsity $S$:} Figure \ref{LolRhoEst} shows
the estimation of the sparsity provided by LOL as a function of the
effective sparsity $S$.
 For small $S$ ($\rho \leq 5\%$), LOL
 is excellent because it estimates exactly (with no error)
the sparsity $S$ for all the studied indeterminacy levels. As the
sparsity $S$ increases, LOL underestimates the parameter $S$. For a
given sparsity value, the underestimation becomes weaker as the
indeterminacy level $\delta$ decreases. Comparing Figure
\ref{LolEY2rho} and Figure \ref{LolRhoEst}, we observe that the
estimation sparsity is obviously linked to the prediction error
which is not a surprise.

\begin{figure}
\makebox{\includegraphics[width=10cm,height=10cm]{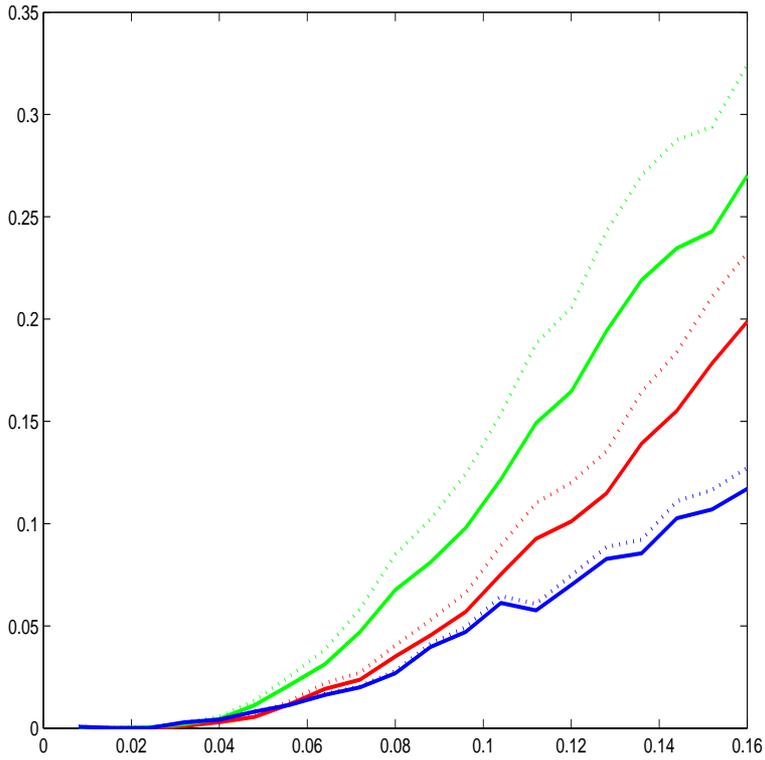}}
\caption{$X-$axis: sparsity rate $\rho$. $Y-$axis: relative
prediction errors for LOL (dot lines) and LOL+ (solid lines).
$\delta=0.4$ (blue color); $\delta=0.75$ (red color); $\delta=0.875$
(green color). The regressors are Gaussian of size $n=250$. $SNR=5$.
\label{LawAdd}}
\end{figure}

{\bf Estimation of the coefficients:} Figure \ref{LawAdd} presents
the improvements provided by LOL$^+$ compared to LOL as a function
of sparsity rate for the prediction error. For all indeterminacy
 and sparsity values, the prediction error decreases using
LOL$^+$ procedure instead of LOL. The improvements are stronger as
both sparsity rate and indeterminacy level increase. The
improvements for the prediction error are observed as $\rho$
increases given all studied indeterminacy levels $\delta$.
Obviously, the estimated sparsity in the same for both procedures
LOL and LOL$^+$.

{\bf Ultra high dimension:}  Table \ref{tabveryhigh} shows the
prediction error for ultra high dimension as $p=5000$, $p=10000$;
$p=20000$ and for two different values of $n=400$ and $n=800$. For
small sparsity levels ($S=5$, $10$, $20$), the performances are
similar even in a very high dimension as $p=20000$. As in the
previous studies in smaller dimension, for higher sparsity levels
($S=40$, $60$), the performances decrease as the sparsity level or
the indeterminacy increases.

%%%%%%%%%%%%%%%%%%%%%%%%%%%%%%%%%%%%%%%%%%%%%%%%%%%%%%%
\subsection{Influence of dependence for gaussian design matrices}
 In the simulations, all the predictors do not have
the same influence because some predictors are directly involved in
the model and some others not. Different type of dependency between
the predictors can also be distinguished: dependency between two
predictors involved (or not involved) in the real underlying model,
and dependencies between two predictors: one involved in the model,
the other not. These  dependencies  have not the same impact on the
results. In order to simulate all possible dependencies, we first
extract a $\Phi_{n,2S}$ sub matrix of $\Phi$ defined by
concatenating vertically the $S$ columns of the predictors included
in the model (and associated with non zero coefficients), and $S$
columns between the $(p-S+1)$ predictors chosen at random not
included in the model. $W_1$ is the associated correlation matrix of
$\Phi_{n,2S}$. A new correlation matrix $W_2$ is then built by
choosing randomly $5\%$ or $20\%$ of the correlations in $W_1$ and
replacing their original value with random values of the form
$(-1)^b u$, where $b$ is drawn from a Bernoulli distribution with
parameter $0.5$ and $u$ from an uniform distribution between $[0.90;
0.95]$ such a way that $W_2$ presents some high correlations between
the $2S$ selected predictors. Since the correlation matrix of the
$2S$ columns of $Z:=\Phi_{n,2S} W_1^{-\frac{1}{2}}
W_1^{\frac{1}{2}}$ is then $W_2$, we replace the previously removed
columns of $\Phi$ by the columns of $Z$.

\begin{figure}
\makebox{\includegraphics[width=10cm,height=10cm]{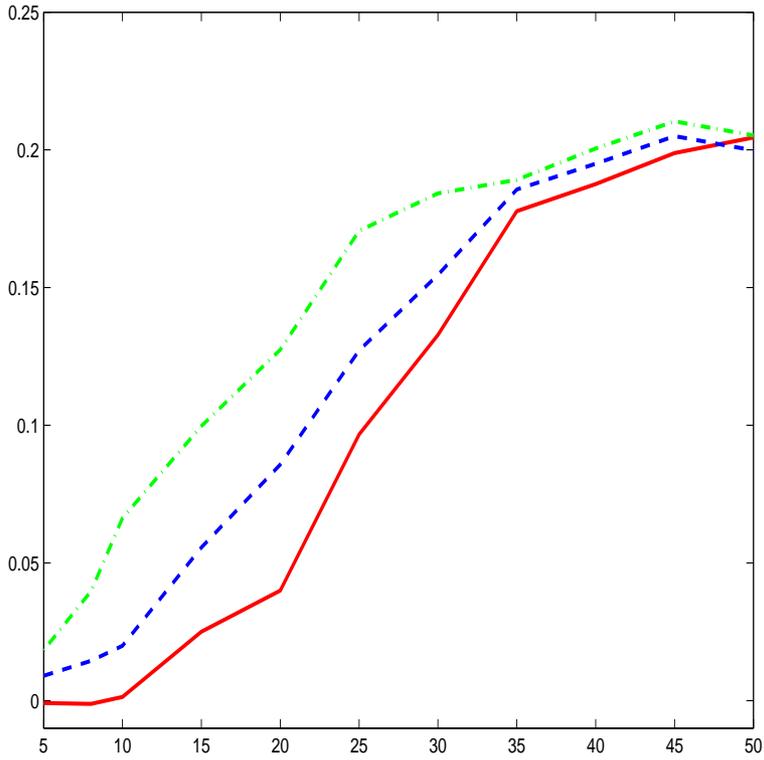}}
\caption{ $X-$axis: sparsity $S$. $Y-$axis: relative prediction
errors, for LOL with independent regressors  (solid line-red) and
dependent regressors ($5\%$ of dependency, dashdot line-blue;
$20\%$, dashed line-blue). $p=1000$, $n=250$, $K=100$. $SNR=5$.
\label{depstudy}}
\end{figure}

Figure \ref{depstudy} compares the prediction error for both
dependent and independent cases. As expected, some dependency
between the predictors damages the performances of LOL. When the
sparsity increases, the impact of dependency seems to play a lower
impact on the prediction error.

\subsection{Impact of the family distribution  of the design matrix}
\begin{figure}
\makebox{\includegraphics[width=10cm,height=10cm]{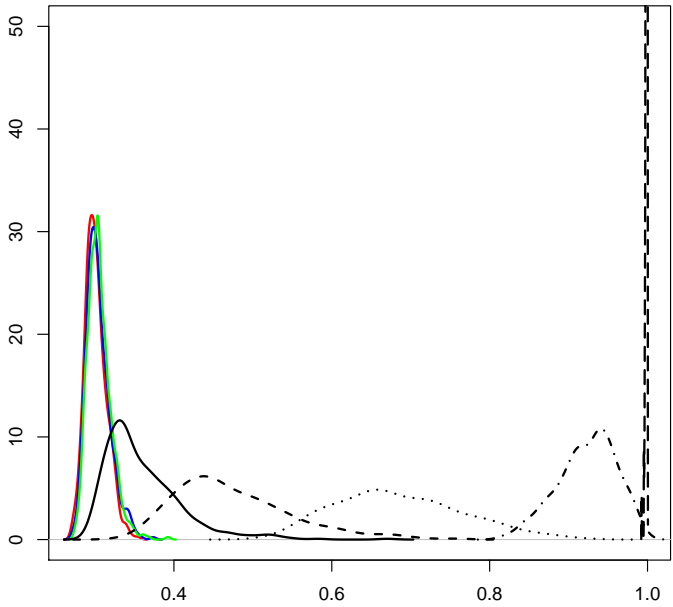}}
\caption{Empirical densities of the coherence $\tau_n$. The
regressors are Gaussien (solid line-red); uniform (solid line-blue);
Bernoulli (solid line-green); Student $5,4,3,2,1$ black lines from
left to right.  $n=250,p=1000$.  \label{Lawcoh}}
\end{figure}

\begin{figure}
\makebox{\includegraphics[width=10cm,height=10cm]{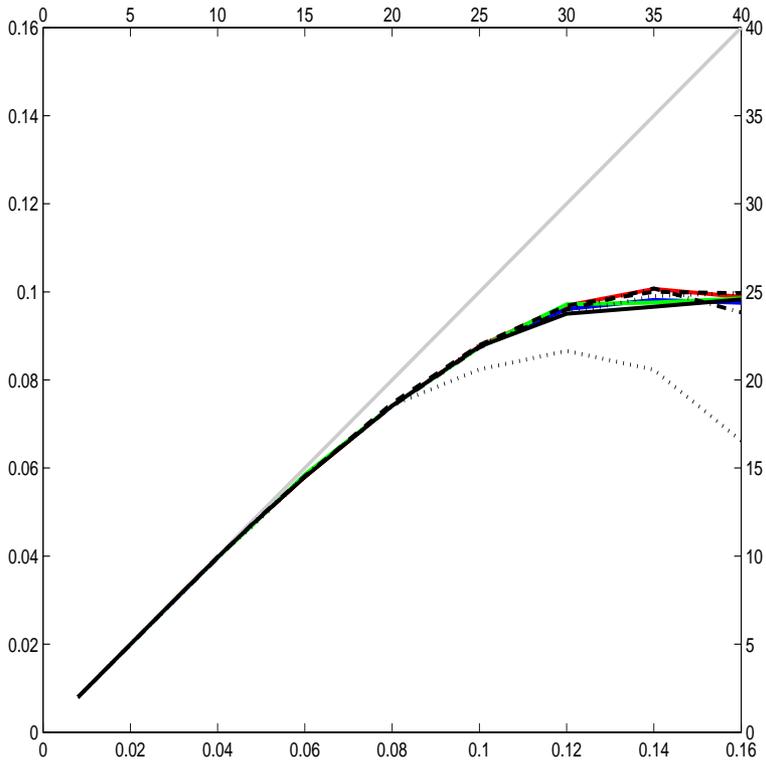}}
\caption{ LOL Sparsity estimation for different distributions for
the predictors. $Gauss$ (solid line-red); $Uniform$ (solid
line-blue); $Bernoulli$: (solid line-green); T(2-5) (black-lines);
T(1) (dot black line). $n=250$, $p=1000$. ($K=200$)
\label{lawstudy}}
\end{figure}

\begin{figure}
\includegraphics[width=10cm,height=10cm]{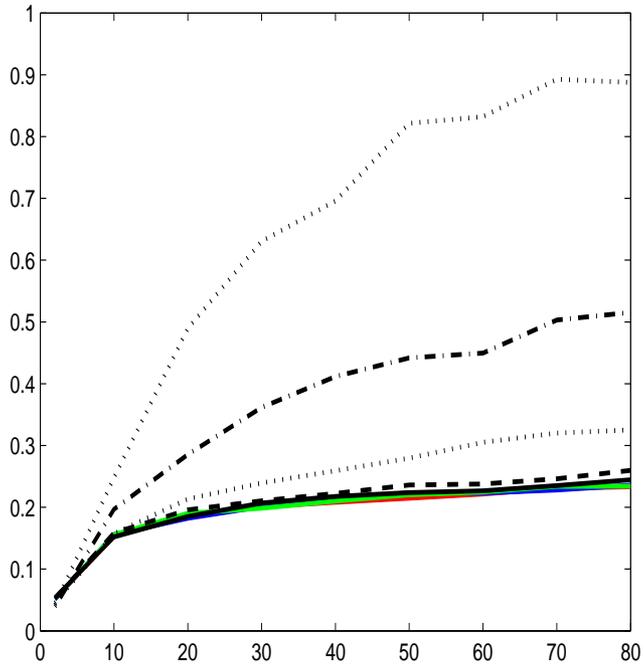}
\caption{$X-$axis: sparsity $S$. $Y-$axis: Coherence $\tau_n$
computed for the $N$ selected Leaders. $Gauss$ (solid line-red);
$Uniform$ (solid line-blue); $Bernoulli$: (solid line-green); T(1)
(dot line-black). $n=250$, $p=1000$. ($K=200$)
\label{leadercohmean}}
\end{figure}

In this section, we investigate the impact of the distribution in a
design matrix with i.i.d. entries. Eight different distributions are
studied: Gaussian ($N(0,1)$), Uniform ($U[-1,1]$), Bernoulli
($B\{-1,+1 \}$) and Student ($T(m)$ with $m \in \{5,4,3,2,1\}$).
Figure \ref{Lawcoh} shows the empirical density of the coherence
$\tau_n$ computed for each law ($n=250$). Similar distributions are
observed for Gaussian, Uniform or Bernoulli laws with a mode of the
coherence equal to $\tau_n=0.30$. For Student's families, a shift of
the mode of the empirical distributions can be observed from left to
right equaled to 0.36 for $T(5)$, 0.47 for T(4), 0.68 for T(3), 0.92
for T(2) to 0.99 for T(1). The prediction errors computed using LOL
are presented in Table $\ref{tablaws}$. For all distributions, the
prediction errors increase with sparsity in average and in
variability. As expected, regarding the coherence value, similar
prediction errors are provided for Gaussian, Uniform, or Bernoulli
laws. For the Student distributions $T(m)$ with parameter $m \geq
2$, the prediction results are also similar to Gaussian
distribution. The Student distribution with $m=1$ shows much higher
prediction errors both in average and variability. Figure
\ref{lawstudy} studies the estimation of sparsity using LOL as a
function of the  sparsity rate $\rho$. All the curves, except the
one for the Student law T(1), are confounded and show similar
behavior as the one observed for gaussian predictors (see Figure
\ref{LolRhoEst} for $\delta=0.25$).  LOL provides similar results
for Gaussian, Uniform, Bernoulli, or Student  laws, $T(m)$ with $m$
large enough. It is amazing to observe that the procedure works fine
even when the empirical coherence $\tau_n$ reaches large values.
However, LOL does not work fine for heavy tailed variables as for
$T(1)$. These results can be explained analyzing Figure
\ref{leadercohmean} which shows the coherence of the matrix
restricted to the $N$ selected leaders. This restricted coherence is
much lower than the coherence computed on all the predictors. For
the Student $T(1)$ law, $\tau_n=0.99$ (see Figure \ref{Lawcoh})
while the coherence restricted to the leaders is $0.3$ (see Figure
\ref{leadercohmean} by instance for $S=10$). LOL  provides also good
results even when the global coherence approaches $1$. It seems then
that the practical results are much more optimistic than the
theoretical ones, although they show deteriorations under high
coherence. Conclusions would be that it could be interesting to find
new measures of collinearity to reflect better the performances of
the method. This is true in general for all the methods concerned
with high dimension.

\subsection{Comparison with other two-step procedures}
In this part, the performances of LOL are compared with the
performances of other two-step procedures which have been
practically studied. The first one referred as SIS-Lasso is coming
from \cite{fanlv}: the selection step called SIS is followed by the
Lasso procedure. The second one called Lasso-Reg, is proposed in
\cite{candesplan}. First, the Lasso algorithm performs the selection
of the leaders and then, the coefficients are estimated with a
regression. For simplicity of the presentation, we do not include
the results provided by greedy algorithms.

 The performances of the three
procedures (LOL, SIS-Lasso, Lasso-Reg) are here studied over a large
range of sparsity  in order to cover previous results already
presented in \cite{fanlv} and \cite{candesplan} for different
sparsity.
  The number of initial predictors is
$p=1000$ and the number of observations $n=200$. This experimental
design allows us to analyze extremely small sparsity values ($10
\leq S \leq 20$) (as in \cite{fanlv}) as well as values as large as
$S=60$ (as in \cite{candesplan}). For the Lasso procedures, the
regularization parameter is chosen by cross validation. Different
signal over noise ratio are studied ($SNR=10$, $5$, $2$).

Table \ref{methodstudy} presents the relative prediction error as
defined for i.i.d. gaussian matrices but similar results are
obtained with uniform, or Bernoulli distribution. Different cases of
signal over noise ratio are studied ($SNR=10, 5, 2$). The
performances of the procedures appear to depend on the sparsity and
on the signal over noise ratio. For small sparsity levels, ($S=10$),
all the procedures perform extremely well and the relative
prediction error is similar to the inverse of the signal over noise
ratio. For middle sparsity levels ($20 \leq S \leq 30$), Lasso-Reg
performs better than the others ones when the signal over noise
ratio is high ($SNR=10$ or $5$). In this case, Lasso-Reg seems to be
more efficient to select (during the first step) the leaders than
both SIS-Lasso and  LOL. For a low signal over noise ratio
($SNR=2$), LOL performs better than Lasso-Reg. The performances of
SIS-Lasso and LOL are globally similar.

For largest values of the sparsity level $S \geq 50$, it appears
that SIS-Lasso and LOL are better than Lasso-Reg for middle values
of the signal over noise ratio.

We conclude that LOL has a special gain over the other procedures
when the SNR is small or when the sparsity $S$ is high.

\subsection{LOL in Boston}
In order to illustrate the performances of LOL on real data, we
revisit the Boston Housing data (available from the {\it UCI machine
learning data base repository}: http://archive.ics.ucfi.edu/ml/) by
fitting predictive models using LOL. The original Boston Housing
data have one continuous target variable $Y$ (the median value of
owner-occupied homes in USD)  and $p_0=13$ predictive variables over
$n=506$ observations which are randomly split into two subsets: one
training set with $75\%$ of observations and one test set with the
remaining $25\%$ observations.

 In view to test our
procedure, we consider the linear regression method as a benchmark
and denote $E^{Reg}$ the prediction error computed on the test set
while the estimated model is computed on the training set.

 The data are 'dived' in a high dimensional space
of size $p=2113$ by adding $300$ independent random variables of
seven different laws: Normal, lognormal, Bernoulli, Uniform,
exponential with parameter 2, Student $T(2)$, $T(1)$ in equal
proportion. This set of laws is chosen to mimic the different
underlaying laws of the 13 original variables. LOL is applied on the
training set and the error of prediction  $E^{LOL}$ is computed on
the test set. This procedure is repeated $K=100$ times using re
sampling, and the prediction errors are then averaged to compute the
performances on the training and test sets. Observe that in this
example, the indeterminacy level and the sparsity rate are quite low
equal to $\delta=0.18$ and $\rho=0.034$. The coherence is quite high
equal to $\tau_n=0.98$.

The results are the following
$$E^{LOL}=0.245 (0.05) \quad \mbox{ and }\quad  E^{Reg}=0.266 (0.04)$$
and LOL appears to work in this case very well because similar
prediction errors are obtained even from a high dimensional space
$p=2113$  as using a regular linear regression in $p_0=13$
dimensions.

%%%%%%%%%%%%%%%%%%%%%%%%%%%%%%%%%%%%%%%%%%%%%%%%%%%%%%
\section{Proofs}\label{section_proofs}
%%%%%%%%%%%%%%%%%%%%%%%%%%%%%%%%%%%%%%%%%%%%%%%%%%%%%%

First, we state preliminary results and next we prove Theorem
\ref{merceradapt} and Theorem \protect{ \ref{lqresults}} as a
consequence of Theorem \protect{\ref{merceradapt}}. The proofs of
the preliminaries are postponed in the appendix.

 For any subset of indices ${\cI}\subset
\{1,\ldots,p\}$, $V_{\cI}$ denotes the subspace of $\R^n$ spanned by
the columns of the extracted matrix $\Phi_{|\cI}$ and $P_{V_{ \cI}}$
denotes the projection over $V_\cI$ (in euclidean sense in $\R^n)$.
Set $\bar\alpha(\cI)$ the vector of $\R^{\#(\cI)}$ such that
$\Phi_{|\cI} \bar \alpha(\cI):=P_{V_\cI} [\Phi\alpha]$. Obviously,
as soon as $\#(\cI)\le N$, we get
$$\bar\alpha(\cI)=(\Phi_{|\cI}^t\Phi_{|\cI })^{-1}\Phi_{|\cI}^t
\Phi\alpha.$$ As well, set $\hat \alpha(\cI)$ such that
$\Phi_{|\cI}\hat \alpha(\cI):=P_{V_{\cI}}( Y)$.

%%%%%%%%%%%%%%%%%%%%%%%%%%%%%%%%%%%%%%%%%%%%%%%%%%%%%%%
\subsection{Preliminaries}
%%%%%%%%%%%%%%%%%%%%%%%%%%%%%%%%%%%%%%%%%%%%%%%%%%%%%%%
The preliminaries contain three essential results for the subsequent
proof. The first proposition describes the algebraic behavior of the
euclidean norm of $\hat \alpha(\cB)-\alpha$ when the vector is
restricted to a (small) set of indices. The second lemma is a
consequence of the RIP property and gives an algebraic equivalent
for the  projection norm of vectors over spaces of small dimensions.
The second proposition (third result) describes the concentration
property for  projections norms of the vector of errors. This
proposition is our major ingredient for proving all the exponential
bounds. Note that it also incorporates the case where the projection
has a possibly random range.

\begin{proposition}\label{projection} Let $\cI$ be a subset
of the leaders indices set $\cB$. Then
\begin{align*}
\sum_{\ell\in \cI}(\widehat{\alpha_\ell} (\cB)-\alpha_\ell)^2&\le
\kappa(\alpha) \#(\cI)\tau_n^2\left( 1+\|P_{V_{\cB}}[\eps]\|_{n}^2
\right) +\frac{3c_0^2}{1-\nu}\frac Sn+
\frac{3}{1-\nu}\|P_{V_{\cI}}[\eps]\|_{n}^2
\end{align*}
where
$$
\kappa(\alpha)=\frac{3(5+\nu)}{(1-\nu)}\|\alpha\|^2_{l^1(p)} \vee
\frac{3}{(1-\nu)}.
$$
\end{proposition}

 %%%%%%%%%%%%%%%%%%%%%%%%%%%%%%%%%%%%%%%%%%%%%%%%%%%%%%%
 \begin{lemma}\label{projobis}
 Let $\cI$ be a subset of $\{1,\ldots,p\}$ satisfying $\#(\cI)\leq N$.
  Then, for any $x \in\R^n$, we get
 \begin{align*}
 (1+\nu)^{-1}\sum_{\ell\,\in \cI}\;\frac 1n\left(\sum_{i=1}^n
  x_i\Phi_{i\ell}\right)^2\le
 \|P_{V_\cI}x\|_{l^2(n)}^2\le (1-\nu)^{-1}
 \sum_{\ell\,\in \cI}\;\frac 1n\left(\sum_{i=1}^n
  x_i\Phi_{i\ell}\right)^2\,.
 \end{align*}
 \end{lemma}

  %%%%%%%%%%%%%%%%%%%%%%%%%%%%%%%%%%%%%%%%%%%%%%%%%%%%%%%
 \begin{proposition}\label{chi2bis}
 Let $\cI$ be a non random subset of $\{1,\ldots,p\}$ such that
  $\#(\cI)\le n_{\cI}$, where $n_{\cI}$ is a deterministic quantity, then
   \begin{align}\label{expon}
 P\left(\frac{1}{\sigma^2}\|P_{V_{\cI}}[\eps]\|_{n}^2\geq
 \mu^2\right)&\leq \exp\left(-n\mu^2/16\right)
 \end{align}
 for any $\mu$ such that $\mu^2\geq 4\,n_{\cI}/n$.
 If now $\cI$ is a random subset of $\{1,\ldots,p\}$ such that
  $\#(\cI)\le n_{\cI}$, where $n_{\cI}$ is a deterministic quantity, then
\eref{expon} is still true but
 for any $\mu$ such that $\mu^2\geq 16\,n_{\cI}\log{p}/n$.
\end{proposition}

Proposition \ref{projection} and Proposition \ref{chi2bis}  are
proved in the appendix as well as
 Lemma \ref{projobis}.

%%%%%%%%%%%%%%%%%%%%%%%%%%%%%%%%%%%%%%%%%%%%%%%%%%%%%%%
\subsection{Proof of Theorem \protect{\ref{merceradapt}}}
%%%%%%%%%%%%%%%%%%%%%%%%%%%%%%%%%%%%%%%%%%%%%%%
For sake of simplicity, and without loss of generality, we assume
that the $N$ largest $\alpha_\ell$'s have their indices in
$\{1,\ldots,N\}$. We have
\begin{align*}d(\hat\alpha^*,\alpha)
&\le \| u\|_{n}+
\|\sum_{\ell=1}^p(\hat\alpha_\ell^*-\alpha_\ell)\Phi_{\bullet\ell}\|_{n}
\;.
\end{align*}
Recall that $\cB$ is the set of the indices of the leaders. Then
\begin{align*}
\|\sum_{\ell=1}^p(\hat\alpha_\ell^*-\alpha_\ell)\Phi_{\bullet\ell}\|_{n
}&\leq \|\sum_{\ell=1}^p\I\{\ell\in
\cB\}(\hat\alpha_\ell^*-\alpha_\ell)\Phi_{\bullet\ell}\|_{n }
+\|\sum_{\ell=1}^p\I\{\ell\not \in
\cB\}\,\alpha _\ell\Phi_{\bullet\ell}\|_{n }\\
&:= I \;(\mbox{In})\;+O \;(\mbox{Out})\;.
\end{align*}
We split $I$ into four terms by observing that :
\begin{align*}
1&= \I \{|\hat\alpha_\ell(\cB)|\ge \lan(2)\} \big\{
\I\{|\alpha_\ell| \ge \lan(2)/2\}+ \I\{|\alpha_\ell|< \lan(2)/2\}
\big\}
\\
&+\I\{|\hat\alpha_\ell(\cB)|< \lan(2)\} \big\{\I\{|\alpha_\ell|\ge
2\lan(2)\}+\I \{|\alpha_\ell|< 2\lan(2)\} \big\}.
\end{align*}
It follows that
\begin{align*}
 I
&\leq \left(\|\sum_{\ell=1}^N\;\I\{\ell\in
\cB\}(\hat\alpha_\ell^*-\alpha_\ell)\Phi_{\bullet\ell}\;\I\{|\alpha_\ell|\ge
\lan(2)/2\}\;\I\{|\hat\alpha_\ell(\cB)|\ge \lan(2)\}\|_{n }\right.\\
& \hspace{2cm}+\left.\|\sum_{\ell=1}^N\;\I\{\ell\in
\cB\}\alpha_\ell\Phi_{\bullet\ell}\;\I\{|\alpha_\ell|\ge
2\lan(2)\}\;\I\{|\hat\alpha_\ell(\cB)|< \lan(2)\}\|_{n }\right)\\
&+\left(
\|\sum_{\ell=1}^p\;\I\{\ell\in\cB\}(\hat\alpha_\ell^*-\alpha_\ell)\Phi_{\bullet\ell}\;\I\{|\alpha_\ell|<
\lan(2)/2\}\;\I\{|\hat\alpha_\ell(\cB)|\ge \lan(2)\}\|_{n }\right.\\
&\hspace{2cm}
 +\left.
\|\sum_{\ell=1}^p\I\{\ell\in\cB\}\alpha_\ell\Phi_{\bullet\ell}\;\I\{|\alpha_\ell|<
2\lan(2)\}\;\I\{|\hat\alpha_\ell(\cB)|\le \lan(2)\}\|_{n }\right)\\
&:=IBB \;(\mbox{InBigBig})\;+ISB
\;(\mbox{InSmallBig})\;+IBS\;(\mbox{InBigSmall})\;+ISS
\;(\mbox{InSmallSmall})\;.
\end{align*}
Note that because  of Assumption (\ref{weak}), the coefficients such
that $|\alpha_\ell|\ge\lan(2)/2$ necessarily have their indices less
than $N$, so some terms in the above sum  have their summation up to
$N$, some others up to $p$. This makes an important difference in
the sequel because Lemma \ref{rho-eucl} can be used in the first
case. Recall the definition of the $\tilde\alpha_\ell$'s given in
Algorithm \ref{algo}
\begin{eqnarray*}
\tilde\alpha_\ell&:=&\frac 1n\,\sum_{i=1}^nY_i\Phi_{i\ell}.
\end{eqnarray*}
We have
\begin{align*}
 O&\leq \|\sum_{\ell=N+1}^p\;\I\{\ell\not \in
\cB\}\,\,\alpha_\ell\Phi_{\bullet\ell}\|_{n
}+\|\sum_{\ell=1}^N\;\I\{\ell\not \in
\cB\}\,\alpha_\ell\Phi_{\bullet\ell}\;\|_{n }
 \\
 &\leq \|\sum_{\ell=N+1}^p\;\I\{\ell\not \in
\cB\}\,\,\alpha_\ell\Phi_{\bullet\ell}\|_{n }+
 \|\sum_{\ell=1}^N\;\I\{\ell\not \in
\cB\}\,\alpha_\ell\Phi_{\bullet\ell}\;\I\{|\alpha_\ell|\le 2\lan(1)\} \|_{n }\\
&+\|\sum_{\ell=1}^N\;\I\{\ell\not \in
\cB\}\,\alpha_\ell\Phi_{\bullet\ell}\;\I\{|\alpha_\ell|\ge
2\lan(1)\}\;\I\{|\tilde\alpha_\ell|\le \lan(1)\}
 \|_{n }\\&+\|\sum_{\ell=1}^N\;\I\{\ell\not \in
\cB\}\,\alpha_\ell\Phi_{\bullet\ell}\;\I\{|\alpha_\ell|\ge
2\lan(1)\}\;\I\{|\tilde\alpha_\ell|\ge \lan(1)\}
 \|_{n }\\
&:=Ob\;(\mbox{OutBias})\;+OS
\;(\mbox{OutSmall})\;+OBS\;(\mbox{OutBigSmall})\;+OBB
\;(\mbox{OutBigBig})\;
\end{align*}
Using  the Assumption (\ref{h}) on the errors, we get
\begin{align*}
 \| u\|_{n }&\le \sup_{i=1,\ldots,n}|u|\leq
c_0\sqrt{\frac Sn}.
\end{align*}
We deduce that for any $\eta$ such that
 \begin{align}\label{eta1}
\eta^2& >2\, c_0^2\frac Sn,
\end{align}
\begin{align*}P\left(d(\hat\alpha^*,\alpha)\geq \eta\right)
&\le P\left(I+O \geq \eta /2\right)\\
&\le
 P\left(IBB+ISB+IBS+ISS \geq \eta /4\right)+ P\left(OS+OBB+OBS+Ob \geq \eta /4\right).
\end{align*}
Our aim is to prove that each probability term is bounded by
$\exp-\gamma n\eta^2$ for any $$\eta^2 \geq D\left(\frac
{S\log{p}}{n}\vee S\tau_n^2\right)$$ where the constants $\gamma$
and $D$ have to be determined. To do this, basically, we
 study each term separately and prove that (up to constants)
either it can be directly bounded, or it reduces to a random term
whose probability of excess can be bounded using Proposition
\ref{chi2bis}.
%%%%%%%%%%%%%%%%%%%%%%%%%%%%%%%%%%%%%
\subsubsection{Study of $IBB$ and $ISB$}
%%%%%%%%%%%%%%%%%%%%%%%%%%%%%%%%%%%%%
Denote by $\cT$ the (non random) set of indices
$\{\ell=1,\ldots,N,\, |\alpha_\ell|\geq \lan(2)/2\}$ which verifies
$\#\cT\leq S$ by Assumption \ref{weak}. Observe that
\begin{eqnarray*}
\left\{ \begin{array}{l}
|(\hat\alpha(\cB))_\ell|<\lan(2)\\
|\alpha_\ell|>2\lan(2)\end{array}\right.&\Longrightarrow&
|(\hat\alpha(\cB))_\ell|<\lan(2)<|\alpha_\ell-(\hat\alpha(\cB))_\ell|
\end{eqnarray*}
Using Lemma \ref{rho-eucl}, we deduce that
\begin{eqnarray*}
 {ISB^2}&\le& (1+\nu)\sum_{\ell\in \cT\cap\cB}|(\hat\alpha(\cB))_\ell+
 (\alpha_\ell-(\hat\alpha(\cB))_\ell)|^2 \;\I\{|(\hat\alpha(\cB))_\ell|\le
 |\alpha_\ell-(\hat\alpha(\cB))_\ell|\}\;\I\{|\alpha_\ell|\ge 2\lan(2)\}\\
&\le&4(1+\nu)\sum_{\ell\in \cT\cap\cB}(\alpha_\ell-
(\hat\alpha(\cB))_\ell )^2.
 \end{eqnarray*}
Using again Lemma \ref{rho-eucl}, it follows that
\begin{eqnarray*}
 {ISB^2}+IBB^2
&\le&5(1+\nu)\sum_{\ell\in \cT\cap\cB}(\alpha_\ell-
(\hat\alpha(\cB))_\ell )^2.
 \end{eqnarray*}
We apply Proposition \ref{projection}
\begin{eqnarray*}
 {ISB^2}+IBB^2
&\le&5(1+\nu)\left(\kappa(\alpha)S\tau_n^2
\left(1+\,\|P_{V_{\cB}}[\eps]\|_{n}^2\right)
+\frac{3c_0^2}{1-\nu}\frac Sn+\frac{3}{1-\nu}
\|P_{V_{\cT}}[\eps]\|_{n}^2\right).
 \end{eqnarray*}
We use now Proposition \ref{chi2bis}: first with the non random set
$\cT$ satisfying $\#(\cT)\leq S$, secondly with the random set $\cB$
such that $\#(\cB)\leq N$. For this second part, we use
 the last part of Proposition \ref{chi2bis}, which yields
 an additional logarithmic factor. We obtain
\begin{align*}
P\left(IBB+ISB \geq \eta /8\right)&\leq
P\left(\frac{1}{\sigma^2}\,\|P_{V_\cT}\eps \|_{n}^2 \geq
\eta^2\,(1-\nu)/( 7680(1+\nu)\sigma^2) \right) \\
&+P\left(\frac{1}{\sigma^2}\,\|P_{V_\cB}\eps \|_{n}^2 \geq
\eta^2(S\tau_n^2)^{-1}/( 2560 (1+\nu)\kappa(\alpha)\sigma^2) \right)
\\
&\leq 2\exp\left\{-n\eta^2(1-\nu)/( 30720(1+\nu)\sigma^2\right)
\end{align*}
since  $S\leq N=\nu/\tau_n$ and as soon as
\begin{align}\label{eta2}
\eta^2 &\geq
2560\left((1+\nu)\kappa(\alpha)S\tau_n^2\vee\frac{c_0^2}{1-\nu}\frac{S}{n}\right)
\vee  30720 \frac{(1-\nu)}{(1+\nu)}\sigma^2\frac Sn\nonumber
\\
&\hspace{1cm}\vee 40960 \nu(1+\nu)\nu\,\kappa(\alpha)\sigma^2
\frac{S\tau_n\,\log{p}}{n}.
\end{align}

%%%%%%%%%%%%%%%%%%%%%%%%%%%%%%%%%%%%%
\subsubsection{Study of $Ob$}

Since the $\Phi_{\bullet\ell}$'s are normalized vectors and because
of the
 definition of the coherence, we get
\begin{eqnarray*}
Ob&\leq&\|\sum_{\ell\geq N+1}\alpha_\ell\Phi_{\bullet\ell}\|_{n }\\
&\leq&\left[\sum_{\ell\geq N+1}\alpha_\ell^2+\tau_n
\left(\sum_{\ell\geq N+1}|\alpha_\ell|\right)^2\right]^{1/2}
\\
&\leq&\left[\sum_{\ell\geq N+1}\alpha_\ell^2\,\I\{|\alpha_\ell|\leq
\lan(2)/2\}+\tau_n \left(\sum_{\ell\geq
N+1}|\alpha_\ell|\right)^2\right]^{1/2}.
\end{eqnarray*}
As $\lan(2)\leq \lan(1)$ and using Assumption (\ref{l2q}) and
Assumption (\ref{l1q}), we obtain
$$Ob\le
\,(c_1+c_2)\sqrt{\frac{S\log{p}}{n}}
$$ which implies that $Ob\leq \eta/16$ as soon as
\begin{align}\label{eta3}
\eta^2&
>256(c_1+c_2)^2\;\frac{S\log{p}}{n}.
\end{align}

%%%%%%%%%%%%%%%%%%%%%%%%%%%%%%%%%%%%%%%
\subsubsection{Study of $ISS$ and $OS$} As $\lan(1)\ge \lan(2)$, using
successively Lemma \ref{rho-eucl} and Assumption \eref{l2q}, we have
\begin{align*}
ISS&\leq
Ob+\|\sum_{\ell=1}^N\I\{\ell\in\cB\}\alpha_\ell\Phi_{\bullet\ell}\;\I\{|\alpha_\ell|\le
2\lan(2)\}\;\I\{|\hat\alpha_\ell(\cB)|\le \lan(2)\}\|_{n }\\
&\leq
Ob+(1+\nu)^{1/2}\,\left(\sum_{\ell=1}^N\alpha_\ell^2\;\I\{|\alpha_\ell|\le
2\lan(1)\}\right)^{1/2}\\
&\leq (c_1+c_2)\sqrt{\frac{S\log{p}}{n}}+
(1+\nu)^{1/2}c_2\sqrt{\frac{S\log{p}}{n}}.
\end{align*}
This implies that $ISS\leq \eta/16$ as soon as
\begin{align}\label{eta4}
\eta^2&>512\,\left((c_1+c_2)^2+(1+\nu)c_2^2\right)\,\frac
{S\log{p}}{n}.
\end{align}
In the same way, $OS\leq \eta/16$.
%%%%%%%%%%%%%%%%%%%%%%%%%%%%%%%%%%%%%%%%%%%%%%%%%%%%%%%%%%%%%%%%%%%%%%%%%%%%%%%%%%%%%%%%%%%%%%%%%%%%%%%%%%%%%%%%%%%%%%%%%%%%%%%%%%%%%%%%%%%%%%%%%%%%%%%%%%%%%%%%%%%%%%%%%%%%%%%%%%%%%%%%%%%%%%%%%%%%%%%%%%%%%%%%%%%%%%%%%%%%%%%%%%%%%%%%%%%%%%%%%%%%%%%%%%%%%%%%%%%%%%%%%%%%%%%%%%%%%%%%%%%%%%%%%%%%%%%%
\subsubsection{Study of $OBS$}
Using the model and the definition of
 $\widetilde{ \alpha_{\ell}}$ given in Algorithm \ref{algo}, we get
$$\widetilde{ \alpha_{\ell}}=
\frac
1n\sum_{i=1}^n[\sum_{m=1}^p\alpha_m\Phi_{im}+u_i+\eps_i]\,\Phi_{i\ell}.
$$
Since $\Phi$ has normalized columns, we can write
$$ \alpha_{\ell}=
\frac 1n\sum_{i=1}^n\alpha_\ell\Phi_{i\ell}\Phi_{i\ell}
$$
which implies that
\begin{eqnarray}
|\widetilde{\alpha_{\ell}}-\alpha_{\ell}| &\leq&\left|\frac
1n\sum_{i=1}^n\left(\sum_{m=1}^p\alpha_m\Phi_{im}\Phi_{i\ell}-\alpha_\ell\Phi_{i\ell}\Phi_{i\ell}\right)\right|
+\left|\frac 1n\sum_{i=1}^nu_{i}\Phi_{i\ell}\right| +\left|\frac
1n\sum_{i=1}^n\eps_i\Phi_{i\ell}\right|\nonumber
\\&\leq&
\left|\sum_{m=1,m\not=\ell}^p\alpha_m\,\frac
1n\sum_{i=1}^n\left(\Phi_{im}\Phi_{i\ell}\right)\right| +\left|\frac
1n\sum_{i=1}^nu_{i}\Phi_{i\ell} \right|
+\left|\frac 1n\sum_{i=1}^n\eps_i\Phi_{i\ell}\right|\nonumber\\
&\leq &\sum_{m=1,m\not=\ell}^p|\alpha_m|\tau_n+\left|\frac
1n\sum_{i=1}^nu_{i}\Phi_{i\ell} \right| +\left|\frac
1n\sum_{i=1}^n\eps_i\Phi_{i\ell}\right|. \label{alpha_alphatilde}
\end{eqnarray}
Recall that $\lan(1)\ge \lan(2)$. We get
$$\I\{|\widetilde{\alpha_{\ell}}|< \lan(1)\}\;\I\{|\alpha_\ell|\ge 2\lan(1)\}\quad\le\quad
\I\{|\alpha_\ell- \widetilde{\alpha_{\ell}}|\ge \lan(1)\ge
|\widetilde{\alpha_{\ell}}|\}\;\I\{|\alpha_\ell|\ge \lan(2)/2\}.
$$
Hence, using Lemma \ref{rho-eucl}, it follows
 \begin{eqnarray}
OBS^2&\le& (1+\nu)\sum_{\ell=1}^N(\alpha_\ell-
\widetilde{\alpha_{\ell}}+\widetilde{\alpha_{\ell}} )^2
\;\I\{|\alpha_\ell- \widetilde{\alpha_{\ell}}|\ge \lan(1)\ge
|\widetilde{\alpha_{\ell}}|\}\,\I\{|\alpha_\ell|\ge
\lan(2)/2\}\nonumber\\
&\le&4(1+\nu)\sum_{\ell=1}^N(\alpha_\ell- \widetilde{\alpha_{\ell}}
)^2\,\I\{|\alpha_\ell|\ge \lan(2)/2\}\label{OBS}.
 \end{eqnarray}
Denote $\cT$ the (non random) set of indices $\{\ell=1,\ldots,N,\;
|\alpha_\ell|\ge \lan(2)/2\}$. Using inequality
(\ref{alpha_alphatilde}), we obtain
 \begin{eqnarray*}
OBS^2&\le&8(1+\nu)\left[\sum_{\ell\in\cT}
(\sum_{\ell^\prime=1}^p|\alpha_{\ell^\prime}|\tau_n)^2+2
\sum_{\ell\in\cT}\left( \left(\frac
1n\sum_{i=1}^nu_i\Phi_{i\ell}\right)^2+\left(\frac
1n\sum_{i=1}^n\eps_i\Phi_{i\ell}\right)^2\right)\right]\\
&:=& OBS_{1}+OBS_{2}+OBS_{3}.
 \end{eqnarray*}
By Assumption (\ref{weak}), we have $\#\cT\leq S$ implying
$$
OBS_{1}\leq 8(1+\nu)\|\alpha\|_{l^1(p)}^2\,\tau_n^2\, S.
$$
Using Lemma \ref{projobis} and Assumption (\ref{h}) on the errors
$u$, we get
$$
OBS_{2}\leq 16(1+\nu)\, \|P_{V_\cT}[u]\|_{n}^2\le 16(1+\nu)
\,\|u\|_{n}^2\le  16(1+\nu)c_0^2 \frac Sn
$$
and
$$
OBS_{3}\leq 16(1+\nu) \|P_{V_\cT}[\eps]\|_{n}^2.$$ Proposition
\ref{chi2bis} ensures that
 \begin{eqnarray*}
P(OBS>\eta/16 ) &\leq& P\left( \frac{1}{\sigma^2}\,
\|P_{V_\cT}[\eps]\|_{n}^2\geq \eta^2/( 8192\,\sigma^2(1+\nu))\right)\\
&\leq&\exp\left(-n\eta^2/(131072 \,\sigma^2(1+\nu))\right)
 \end{eqnarray*}
as soon as
 \begin{align}
 \label{eta5}
\eta^2&\geq  8192(1+\nu)\;\left(\|\alpha\|_{l^1(p)}^2\,S\tau_n^2\vee
(2c_0^2\vee4\sigma^2)\frac {S}{n}\right)\,.
 \end{align}
 %%%%%%%%%%%%%%%%%%%%%%%%%%%%%%%%%%%%
 \subsubsection{Study of $OBB$}
Observe that the (random) set of indices
$$\cT=\{\ell\not\in \cB,\;|\alpha_\ell|\ge 2\lan(1)\;, |\tilde \alpha_\ell|\ge
\lan(1)\}$$ has no more than $S$ elements (using Assumption
(\ref{weak}) with $\lan(1)\geq \lan(2)$) and is equal to
$\cT_1\cup\cT_2$ where $\cT_1=\cT\cap \{\ell,|\tilde
\alpha_\ell|\le| \alpha_\ell|/2\} $ and $\cT_2=\cT\cap
\{\ell,|\tilde \alpha_\ell|\ge| \alpha_\ell|/2\} $. On the one hand,
we obviously have
\begin{align}\label{T1}\cT_1&\subset \{\ell\not\in \cB,\;|\alpha_\ell|\ge 2\lan(1)\;, |\alpha_\ell|\le 2|\tilde \alpha_\ell- \alpha_\ell|\}.
\end{align} On the other hand,
since $\ell\not\in \cB$ while $|\tilde \alpha_{\ell}|\geq \lan(1)$,
there exist at least $N$ (leader) indices $\ell^\prime$ in
$\{1,\ldots,p\}$ such that $|\tilde \alpha_{\ell^\prime}|\geq
|\tilde \alpha_{\ell}|$. Moreover Assumption \eref{weak} ensures
that there is no more than $S$ indices $\ell^\prime$ such that $|
\alpha_{\ell^\prime}|\ge \lan(1)/2$. Thus, using the fact that
$S<N$, we deduce that there exists at least one index depending on
$\ell$ called $\ell^*(\ell)$ such that
$$
 | \alpha_{\ell^*(\ell)}|\le \lan(1)/2 \mbox{
and } |\tilde \alpha_{\ell^*(\ell)}|\ge |\tilde \alpha_\ell| .
$$
Since $\ell\in \cT_2$, this implies that
\begin{align}\label{T2}
|\tilde \alpha_{\ell^*(\ell)}- \alpha_{\ell^*(\ell)}|&\ge |
\alpha_\ell|/4.
\end{align}
 Using (\ref{T1}) and (\ref{T2}), it follows that
\begin{align*}
 OBB^2
 &\le (1+\nu)\left[ \sum_{\ell\in\cT_1} | \alpha_\ell|^2
 \;+\;\sum_{\ell\in\cT_2}| \alpha_{\ell}|^2\right]\\
 &\le (1+\nu)\left[ \sum_{\ell=1}^ N4\, |\tilde \alpha_\ell- \alpha_\ell|^2\;\I\{|\alpha_\ell|\ge 2\lan(1)\}
 \quad+\sum_{\ell\in \; \cT_2} 16\,|\tilde \alpha_{\ell^*(\ell)}-  \alpha_{\ell^*(\ell)}|^2\right]\\
 &:=OBB_1+OBB_2.
 \end{align*} Since $\lan(1)\geq \lan(2)$, $OBB_1$ can be bounded as $OBS$. The computations are exactly the same for the term $OBB_2$ except that the set $\cT_2$ is now random and the conditions on $\eta$ become
 \begin{align}\label{eta7}
\eta^2&\geq   32768
(1+\nu)\;\left(\|\alpha\|_{l^1(p)}^2\,S\tau_n^2\vee 2c_0^2\frac
{S}{n}\vee 16\sigma^2\frac {S\,\log{p}}{n}\right)\,.
 \end{align}
For such an $\eta$, we obtain
 \begin{eqnarray*}
P(OBB>\eta/16 ) &\leq&\exp\left(-n\eta^2/( 524288
 \,\sigma^2(1+\nu))\right).
 \end{eqnarray*}

%%%%%%%%%%%%%%%%%%%%%%%%%%%%%%%%%%%%%%%
\subsubsection{Study of $IBS$}
Note here that the major difficulty lies in the fact that the
summation is not on the set of indices $\ell\le N$ as for the other
terms. Let $\cT,\cT^\prime$ be the subsets of $\{1,\ldots,p\}$
defined as follows
$$\cT=\{\ell\in \cB,\; |\widetilde{\alpha_\ell}|\ge \lan(1),\;
|\hat\alpha_\ell(\cB)|\ge \lan(2),\;|\alpha_\ell|\le \lan(2)/2\}$$
and
$$\cT^\prime=\{\ell\in \cB ,\; |\widetilde{\alpha_\ell}-\alpha_\ell|\ge \lan(1)/2,\;
|\hat\alpha_\ell(\cB)-\alpha_\ell|\ge \lan(2)/2,\;|\alpha_\ell|\le
\lan(2)/2\}$$ and observe that $\cT\subset \cT^\prime$ (using again
that $\lambda_n(2)<\lambda_n(1)$). Denote
$$
K(\cT)=\# \left(\cT\cap\{\ell,\,|\tilde\alpha_\ell-\alpha_\ell|\ge
\lan(1)/2\}\right)
$$
and put $k_0=\lfloor\frac{1}{\textcolor{blue}{8192}
(1+\nu)\kappa(\textcolor{blue}{\alpha})\sigma^2}
\frac{\eta^2}{\lan^2(1)}\rfloor\wedge N$ $k_0=\lfloor \frac{cte
\eta^2}{\lan^2(1)}\rfloor\wedge N$.
 We get
$$
P(IBS>\eta/16)\leq P(IBS>\eta/16 \mbox{ and }K(\cT)\leq
k_0)+P(IBS>\eta/16 \mbox{ and }K(\cT)> k_0):=p_1+p_2.
$$
Notice that $p_2=0$ when $k_0=N$ since $\cT\subset \cB$. To bound
$p_1$, we proceed rather roughly. By Proposition \ref{projection},
we get, for any $k\leq k_0$
\begin{align*}
P(IBS>\eta/16& \mbox{ and }K(\cT)=k)\le
 P\left((1+\nu)\sum_{\ell\in
\cT\cap\cB}(\alpha_\ell-\hat \alpha(\cB)_{\ell})^2\ge \eta^2/256
 \mbox{ and }K(\cT)=k \right)\\
&\leq P\left(\frac {1}{\sigma^2}\,\|P_{V_\cT}\eps\|^2_{n}\ge
\eta^2(1-\nu)/( 15360 (1+\nu)\sigma^2)\mbox{ and
}K(\cT)=k\right)\\&+ P\left(\frac
{1}{\sigma^2}\,\|P_{V_\cB}\eps\|^2_{n}\ge \eta^2(k\tau_n^2)^{-1}/(
5120
(1+\nu)\kappa(\alpha)\sigma^2)\mbox{ and }K(\cT)=k)\right)\\
&\leq \exp\left(-n\eta^2(1-\nu)/(  245760
 (1+\nu)\sigma^2)\right)\\& +
\exp\left(-n\eta^2(k\tau_n^2)^{-1}/( 8192
(1+\nu)\kappa(\alpha)\sigma^2)\right)\\
&\leq 2\exp\left(-n\eta^2(1-\nu)/(  245760
 (1+\nu)\sigma^2)\right)
\end{align*}
because $k\tau_n^2\leq k_0\tau_n^2\leq N\tau_n^2\leq 1$. The
previous bound is valid for any $k\leq k_0$ as soon as
\begin{align*}
\eta^2&\geq 5120
\left(\kappa(\alpha)S\tau_n^2\vee\frac{c_0^2}{1-\nu}\frac
Sn\right)\\ &\hspace{1cm} \vee   61440
 \frac{1+\nu}{1-\nu}\sigma^2
\frac{k_0}{n}\vee 81920(1+\nu)\sigma^2\kappa(\alpha)
\left(\frac{k_0\tau_n^2N\log{p}}{n} \right)
\end{align*}
 which is equivalent to
\begin{align}\label{eta7}
\eta^2\geq
5120\left(\kappa(\alpha)S\tau_n^2+\frac{c_0^2}{1-\nu}\frac Sn\right)
\end{align}
if
\begin{align}\label{la4}
\lan^2(1)\geq cte 6144\frac{1+\nu}{1-\nu}\sigma^2\frac 1n \vee cte
81920\nu(1+\nu)\sigma^2\kappa(\alpha)\frac{\tau_n\log{p}}{n}.
\end{align}
Finally, we get
\begin{align*}
p_1&\le\sum_{k\le k_0}\;\sum_{\cT,K(\cT)=k}P(IBS>\eta/16 \mbox{ and
}K(\cT)=k)\\
&\leq \sum_{k\le k_0}\left(p^k 2\exp\left(-n\eta^2(1-\nu)/(  245760
 (1+\nu)\sigma^2)\right)\,\right)\\
&\leq 2\exp\bigg\{-\left(n\eta^2(1-\nu)/(  245760
 (1+\nu)\sigma^2)\right)\left(1-\frac{k_0\log{p}}{n\eta^2(1-\nu)/(  245760
 (1+\nu)\sigma^2)}\right)\bigg\} \\
&\leq 2\exp\left(-n\eta^2(1-\nu)/(   491520
 (1+\nu)\sigma^2))\right)
\end{align*}
thanks to the choice of $k_0$ and as soon as
\begin{align}\label{la5}
\lan^2(1)\geq
2*245760(cte)^{-1}\frac{1+\nu}{1-\nu}\sigma^2\frac{\log{p}}n.
\end{align}
 To bound $p_2$ (only in the case where $k_0\le N$), we
proceed as above, considering all the (non random) possible sets for
$\cT$. The inclusion $\cT\subset \cT^\prime$
 ensures that
\begin{align*}
p_2&\le\sum_{k\ge k_0}\;\sum_{\cT,K(\cT)=k} P( \sum_{\ell\in
\cT}|\alpha_\ell-\tilde \alpha_\ell|^2\ge k(\lan(1)/2)^2).
\end{align*}
We already have seen that
 \begin{eqnarray*}
\sum_{\ell\in \cT}(\alpha_\ell- \widetilde{\alpha_{\ell}} )^2
&\le&2\tau_n^2\sum_{\ell\in \cT}[\;\sum_{m=1}^p|\alpha_m|\;]^2+4
\sum_{\ell\in \cT}[\frac 1n\sum_{i=1}^nu_i\Phi_{i\ell}]^2+4
\sum_{\ell\in \cT}[\frac 1n\sum_{i=1}^n\eps_i\Phi_{i\ell}]^2
 \end{eqnarray*}
with
 \begin{eqnarray*}
\sum_{\ell\in \cT}[\;\sum_{m=1}^p|\alpha_m|\;]^2\leq
\#(\cT)\,\|\alpha\|_{l^1(p)}^2\quad \mbox{ and }\quad \sum_{\ell\in
\cT}[\frac 1n\sum_{i=1}^nu_i\Phi_{i\ell}]^2\leq (1+\nu)c_0^2\frac Sn
.
 \end{eqnarray*}
It follows
\begin{align*}
p_2 &\le\sum_{1+ k_0\le k\le N}p^k \;P (4\,\sum_{\ell\in\cT} |\frac
1n\sum_{i=1}^n\eps_i\Phi_{il}|^2]\ge k\lan(1)^2/8)
\end{align*}
as soon as
\begin{align}\label{la2}
2\tau_n^2\|\alpha\|_{l^1(p)}^2\leq \lan(1)^2/16\quad \mbox{ and
}\quad 4(1+\nu)c_0^2\frac Sn\leq k_0\lan(1)^2/16.
\end{align}
Recall that $k_0=\lfloor \frac{cte \eta^2}{\lan^2(1)}\rfloor\wedge
N$. Then the second condition is satisfied as soon as
\begin{align}\label{eta8}
\eta^2 &>64 (cte)^{-1}(1+\nu)c_0^2\frac Sn.
\end{align}
Using again Lemma \ref{projobis}, it follows
\begin{align*}
p_2&\le\sum_{1+ k_0\le k\le N}p^k \;P
(\frac{1}{\sigma^2}\;\|P_{V_{\cT}}[\eps]\|_{n}^2\ge k\lan(1)^2/(32(1+\nu)\sigma^2))\\
&\le\sum_{1+ k_0\le k\le N}p^k\; \exp
\left(-n(1+\nu)k\lan(1)^2/(512(1+\nu)\sigma^2)\right)\\
&\le\sum_{1+ k_0\le k\le N}\; \exp
\left(-\left[nk\lan(1)^2/(1024(1+\nu)\sigma^2)\right]
\left[1-\frac{k\log{p}}{nk\lan(1)^2/(1024(1+\nu)\sigma^2)}\right]
\right)
\\
&\le\sum_{1+ k_0\le k\le N}\; \exp
\left(-nk\lan(1)^2/(2048(1+\nu)\sigma^2) \right)
\end{align*}
for
\begin{align}\label{la3}
\lan^2(1)\geq
 2048\sigma^2(1+\nu)\;\frac{\log{p}}{n}.
\end{align}
It follows that
\begin{align*}
p_2 &\le \exp \left(-nk_0\lan(1)^2/(2048(1+\nu)\sigma^2) \right)
\end{align*}
and replacing $k_0$, we conclude that
\begin{align*}
p_2&\leq \exp \left(-\,n\eta^2\,cte/( 2048
\sigma^2(1+\nu)^2)\right).
\end{align*}
%%%%%%%%%%%%%%%%%%%%%%%%%%%%%%%%
\subsubsection{End of the proof}
%%%%%%%%%%%%%%%%%%%%%%%%%%%%%%%%
We now use Assumption \eref{l1} ensuring that $M$ is the radius of
the $l^1-$ ball of the $\alpha$'s to bound $\kappa_1(\alpha)$ by
$(12M^2+5c_0^2)/(1-\nu)^2$. Collecting the conditions \eref{la4},
\eref{la5}, \eref{la2} and \eref{la3} and  on the level $\lan(1)$,
we obtain the constraint
$$
\lan^2(1)\geq \sigma^2(1+\nu)\left(2048\vee \frac{ 491520
}{cte(1-\nu)}\right)\,\frac{\log{p}}{n}\vee32M^2\tau_n^2.
$$
Moreover $\eta$ has to satisfy successively the conditions
(\ref{eta2}), (\ref{eta1}), (\ref{eta3}), (\ref{eta4}),
(\ref{eta5}), (\ref{eta7}) and (\ref{eta8}) leading to the final
condition
$$
\eta\geq D\left(\frac {S\log{p}}{n}\vee S \tau_n^2 \right)
$$
for (revoir)
$$D= \frac{5120}{(1-\nu)^2}(12M^2\vee5c_0^2)
\vee   163840\frac{\sigma^2}{(1-\nu)^2}\vee512(c_1+c_2)^2\vee 65536
(M^2+16\sigma^2).$$ For such an $\eta$, we have
\begin{align*}P\left(d(\hat\alpha^*,\alpha)\geq \eta\right)
&\le  P\left(IBB+ISB+IBS+ISS \geq \eta /4\right)+ P\left(OS+OBB+OBS+Ob \geq \eta /4\right)\\
&\le  P\left(IBB+ISB \geq \eta /8\right)+P\left(IBS \geq \eta
/8-ISS\right)\\&+
 P\left(OBB \geq \eta /8-OS\right)+ P\left(OBS \geq \eta /8-Ob\right)
\end{align*}
which is bounded by $8\exp(- n\eta^2/\gamma)$ for
$$
\gamma=C(1+\nu)\sigma^2(1+(1+\nu)\sigma^2)
$$
where $C$ is an universal numerical constant.

%%%%%%%%%%%%%%%%%%%%%%%%%%%%%%%%%%%%%%%%%%%%%%%%%%%%%%
\subsection{Proof of Theorem \protect{\ref{lqresults}}}
%%%%%%%%%%%%%%%%%%%%%%%%%%%%%%%%%%%%%%%%%%%%%%%%%%%%%%
To prove that Theorem \ref{lqresults} is a consequence of Theorem
\ref{merceradapt}, we need  to prove
$$
B_q(M)\subset V\left(M^q(T_3/2)^{-q}\,\left(\frac n {\log
p}\right)^{q/2},M\right)\quad\mbox{ and }\quad B_0(S,M)\subset
V(S,M)
$$
for $(c_0,c_1,c_2)$ to be specified. First, assume that $\alpha\in
B_q(M) $ for $q\in(0,1]$. Since $q\leq 1$, we have
$\|\alpha\|_{l^1(p)}\leq \|\alpha\|_{l^q(p)}<M$ and \eref{l1} is
satisfied. Since $\lan(2)\geq T_3\sqrt{\log{p}/n}$ and using Markov
Inequality, we get
\begin{align*}\#\left\{\ell=1,\ldots,p,\;|\alpha_\ell|\ge \lan(2)/2\right\}&\le \#\left\{\ell=1,\ldots,p,\;|\alpha_\ell|\ge  \frac
{T_3}2\sqrt{\frac {\log p}n}\right\}\\&\le M^q\left( \frac
{T_3}2\sqrt{\frac {\log p}n}\right)^{-q}.
\end{align*}
This proves \eref{weak} with $S=M^q\left(\frac {T_3}2\sqrt{\frac
{\log p}n}\right)^{-q} $. When $q=1$,  assuming that the coherence
$\tau_n$ satisfies $\tau_n\leq c(\log{p}/n)^{1/2}$, observe that
 $$
\sqrt{\frac{S\log
p}{n\tau_n}}=\left(\frac{2M}{T_3}\right)^{1/2}\left(\frac{\log{p}}{n\tau_n^2}\right)^{1/4}
\geq  \left(\frac{2}{cT_3M}\right)^{1/2} M\geq
\left(\frac{2}{cT_3M}\right)^{1/2}\sum_{\ell \ge N}|\alpha_{(\ell)}|
$$  and thus \eref{l1q} is verified for $c_1=(cMT_3/2)^{1/2}$. When $q\in(0,1)$,
using again Markov inequality, we get
$$\ell = \#\left\{\ell=1,\ldots,p,\;|\alpha_\ell|\ge
|\alpha_{(\ell)}|\}\le  M^q\;|\alpha_{(\ell)}|^{-q}\right\}
$$ leading to  the bound
$|\alpha_{(\ell)}|\le M {\ell^{-1/q}}$. Recall that $N\tau_n=\nu$.
Thus, for $q\in\; (0,1)$
\begin{align*}\sum_{\ell \ge N}|\alpha_{(\ell)}|
& \le \sum_{\ell \ge N}\ell^{-1/q} M \le M N^{1-1/q}\\
&= M \nu^{1-1/q}\;\tau_n^{1/q -1}\leq M \nu^{1-1/q}\;
\left(\frac{n\,\tau_n^{2/q -1}}{S\log
p}\right)^{1/2}\;\sqrt{\frac{S\log p}{n\tau_n}}.
\end{align*}
Notice that
\begin{align*}
\frac{n\,\tau_n^{2/q-1}}{S}&\leq
M^{-q}(T_3/2)^q\frac n{\log p}\tau_n^{2/q-1}\left(\frac{\log{p}}{n}\right)^{q/2}\\
&\leq \left(c^{2/q-1}M^{-q}(T_3/2)^q\right)\,\left(\frac{\log
p}n\right)^{(1-q)(2-q)/2q}.
\end{align*}
 is bounded by a constant when $\log p/n\le c^\prime$.
 This implies \eref{l1q}.
 Now, we get
\begin{align*}
\sum_{\ell=1}^p|\alpha_{\ell}|^2I\{|\alpha_{\ell}|\le 2\lan(1)\}&\le M^q\left(2\lan(1)\right)^{2-q}\le M^q\left(2T_4\sqrt{\frac{\log p}n}\right)^{2-q}\\
&\le
\frac{(2T_4)^{2-q}}{(T_3/2)^{-q}}\;\left(\left(M^q(T_3/2)^{-q}\right)
\left(\frac{\log{p}}{n}\right)^{-q/2}\right)\;\frac{\log p}n \le
c_2^2\;S\frac{\log p}n
\end{align*}
 which proves \eref{l2q} with $c_2^2=(2T_4)^{2-q}(T_3/2)^{q}$. This ends the proof of Theorem \ref{lqresults} when $q\in (0,1]$.
We finish with the case where $\alpha$ belongs to $B_0(S,M)$ which
is very  simple since we have \eref{l1} and \eref{weak} for free.
\eref{l1q} is obviously true with $c_1=0$ and \eref{l2q} is true
with $c_2=4T_4$ because there are only $S$ non zero coefficients and
thus
$$\sum_{\ell =1}^p|\alpha_{\ell}|^2\,\I\{|\alpha_{\ell}|\le
2\lan(1)\}\le S (2\lan(1))^2\le 4T_4\frac{S\log p}n.$$

%%%%%%%%%%%%%%%%%%%%%%%%%%%%%%%%%%%%%%%%%%%%%%%%%%%%%%%
\section{Appendix}
%%%%%%%%%%%%%%%%%%%%%%%%%%%%%%%%%%%%%%%%%%%%%%%%%%%%%%%

 Recall that
$\bar\alpha(\cI)$ is the vector of $\R^{\#(\cI)}$ such that
$\Phi_{|\cI} \bar \alpha(\cI):=P_{V_\cI} [\Phi\alpha]$. As soon as
$\#(\cI)\le N$,
$$\bar\alpha(\cI)=(\Phi_{|\cI}^t\Phi_{|\cI })^{-1}\Phi_{|\cI}^t
\Phi\alpha.$$ As well,  $\hat \alpha(\cI)$ has been defined by
$\Phi_{|\cI}\hat \alpha(\cI):=P_{V_{\cI}}( Y)$. Using the setting
 (\ref{model}), we get
\begin{align}\label{fifi}
\Phi_{|\cI}\hat \alpha(\cI)
 =P_{V_{\cI}}[\Phi\alpha+u+\eps]=\Phi_{|\cI} \bar
\alpha(\cI)+P_{V_{\cI}}[u+\eps].
\end{align}
%%%%%%%%%%%%%%%%%%%%%%%%%%%%%%%%%%%%%%%%%%%%%%%%%%%%%%%
\subsection{Proof of Lemma \ref{projobis}}
%%%%%%%%%%%%%%%%%%%%%%%%%%%%%%%%%%%%%%%%%%%%%%%%%%%%%%%
Recall that the Gram matrix is defined by
$M(\cI)=n^{-1}\Phi_{|\cI}^t\Phi_{|\cI}$. Let $x\in \R^n$. Since
$$P_{V_\cI}x= \Phi_{|\cI}( n\, M(\cI))^{-1} \Phi_{|\cI}^t x,$$
we obtain
\begin{align*}
\|P_{V_{\cI}}x\|_{l^2(n)}^2 &=( \Phi_{|\cI}^t x)^t\;(n\,
M(\cI))^{-1}\;( \Phi_{|\cI}^t x).
\end{align*}
Applying the RIP Property \eref{cond} and observing that
\begin{align*}
\| \Phi_{|\cI}^t x\|_{l_2(\#(\cI))}^2=( \Phi_{|\cI}^t x)^t\,(
\Phi_{|\cI}^t x)= \sum_{\ell\,\in
\cI}\left(\sum_{i=1}^nx_i\Phi_{\ell i}\right)^2,
\end{align*}  we obtain the announced result.

\subsection{Proof of Proposition \ref{projection}}
We have
\begin{align*}
\sum_{\ell\in \cI}(\alpha_\ell-\widehat{\alpha(\cB)}_\ell)^2&
=\|\alpha_{|\cI}-\widehat{\alpha(\cB)}_{|\cI}\|^2_{l^2(\#(\cI))}\\
&\leq  3\left(\|\alpha_{|\cI}-\bar\alpha(\cI)\|^2_{l^2(\#(\cI))}+
\|\bar\alpha(\cI)-\widehat{\alpha(\cI)}\|^2_{l^2(\#(\cI))} +
\|\widehat{\alpha(\cI)}-\widehat{\alpha(\cB)}_{|\cI}\|^2_{l^2(\#(\cI))}\right)\\
&:=3\left(t_1(\cI)+t_2(\cI)+t_3\right).
\end{align*}
Since
$$
\bar\alpha(\cI)=\alpha_{|\cI}+(\Phi_{|\cI}^t\Phi_{|\cI})^{-1}\Phi_{|\cI}^t\;\Phi_{|\cI^c}\alpha_{|\cI^c}
$$
we get, using twice the RIP Property
\begin{align*}
t_1(\cI)&\leq \frac{1}{1-\nu}\;(\bar\alpha(\cI)-\alpha_{|\cI})^t
\,(n^{-1}\Phi_{|\cI}^t\Phi_{|\cI})\,(\bar\alpha(\cI)-\alpha_{|\cI})\\
&=\frac{1}{1-\nu}\;\frac 1{n^2}\;(\alpha_{\cI^c}^t\Phi_{|\cI^c}^t)
\,\Phi_{|\cI}\,(n^{-1}\Phi_{|\cI}^t\Phi_{|\cI})^{-1}\,\Phi_{|\cI}^t\,(\Phi_{|\cI^c}
\alpha_{|\cI^c})\\
&\leq \frac{1+\nu}{1-\nu}\;\frac
1{n^2}\;\|\Phi_{|\cI}^t\,\Phi_{|\cI^c}\alpha_{|\cI^c}
\|_{l^2(\#(\cT))}^2\\
&= \ \frac{1+\nu}{1-\nu}\;\frac 1{n^2}\sum_{\ell\in
\cI}\left(\sum_{\ell^\prime\in \cI^c} \sum_{i=1}^n
\Phi_{i\ell}\Phi_{i\ell^\prime}\alpha_{\ell^\prime}\right)^2\\
&\leq  \frac{1+\nu}{1-\nu}\; \max_{\ell\not =\ell^\prime}\left|\frac
1{n}\sum_{i=1}^n
\Phi_{i\ell}\Phi_{i\ell^\prime}\right|^2\;\sum_{\ell\in\cI}\left(\sum_{\ell^\prime\in
\cI^c}
| \alpha_{\ell^\prime}|\right)^2\\
&\leq \frac{1+\nu}{1-\nu}\; \#(\cI)\,\tau_n^2 \|
\alpha\|_{l^1(p)}^2.
\end{align*}
  Using Lemma \ref{projobis}, Equality \eref{fifi}, we get
\begin{align*}
t_2(\cI)
 &\leq
\frac{1}{1-\nu}\;\|\Phi_{|\cI}\bar\alpha(\cI)-\Phi_{|\cI}\widehat{\alpha(\cI)}\|^2_{n}\\
&\leq \frac{1}{1-\nu}\;\|P_{V_{\cI}}[\eps+u]\|_{n}^2.
\end{align*}
By Assumption
  (\ref{h}) on the errors $u$, we deduce
\begin{align*}
t_2(\cI)&\leq
\frac{1}{1-\nu}\;\left(\|P_{V_{\cI}}[\eps]\|_{n}^2+c_0^2\frac
Sn\right).
\end{align*}
Now, since $\cI\subset \cB$, we obtain
\begin{align*}
\Phi_{|\cI}\widehat{\alpha(\cI)}-
\Phi_{|\cI}\widehat{\alpha(\cB)}_{|\cI}&=P_{V_\cI}[\Phi_{|\cI}\widehat{\alpha(\cI)}-
\Phi_{|\cI}\widehat{\alpha(\cB)}_{|\cI}]\\
&=P_{V_\cI}[\Phi_{|\cI}\widehat{\alpha(\cI)}-
\Phi_{|\cB}\widehat{\alpha(\cB)}+\Phi_{|\cB\setminus\cI}\widehat{\alpha(\cB)}_{|\cB\setminus\cI}]\\
&=P_{V_\cI}[P_{V_\cI}[\Phi\alpha+u+\eps]-
P_{V_\cB}[\Phi\alpha+u+\eps]
+\Phi_{|\cB\setminus\cI}\widehat{\alpha(\cB)}_{|\cB\setminus\cI}]\\
&=P_{V_\cI}[\Phi_{|\cB\setminus\cI}\widehat{\alpha(\cB)}_{|\cB\setminus\cI}].
\end{align*}
Combining with Lemma \ref{rho-eucl} and Lemma \ref{projobis}, it
leads to
\begin{align*}
t_3&\leq \frac{1}{1-\nu}\;\|\Phi_{|\cI}\widehat{\alpha(\cI)}-
\Phi_{|\cI}\widehat{\alpha(\cB)}_{|\cI}\|^2_{n}\\
&= \frac{1}{1-\nu}\;\|
P_{V_\cI}[\Phi_{|\cB\setminus\cI}\widehat{\alpha(\cB)}_{|\cB\setminus\cI}]
\|^2_{n}\\
&\leq \frac{1}{1-\nu}\;\sum_{\ell\in \cI}\frac 1{n^2}
\left(\sum_{i=1}^n\left(\sum_{\ell^\prime\in
\cB\setminus\cI}\widehat{\alpha(\cB)}_{\ell^\prime}\Phi_{i\ell^\prime}\right)
\Phi_{i\ell}\right)^2\\
&\leq \frac{1}{1-\nu}\;\#(\cI)\;\tau_n^2\,\|\widehat{\alpha(\cB)}\|_{l^1(\#(\cB))}^2\\
&\leq
\frac{4}{1-\nu}\;\#(\cI)\;\tau_n^2\,\left(\|\widehat{\alpha(\cB)}-\bar\alpha(\cB)\|_{l^1(\#(\cB))}^2
+\|\bar\alpha(\cB)-\alpha\|_{l^1(\#(\cB))}^2+\|\alpha\|_{l^1(\#(\cB))}^2\right).
\end{align*}
Now, since $\#(\cB)\leq N$ and $N=\nu/\tau_n$, we obtain
\begin{align*}
t_3 &\leq \frac{4}{1-\nu}\;\#(\cI)\;\tau_n^2\,\left(t_1(\cB)
+t_2(\cB)+\|\alpha\|_{l^1(\#(\cB))}^2\right)\\
&\leq \frac{4}{1-\nu}\;\#(\cI)\;\tau_n^2\,
\left(\frac{1+\nu}{1-\nu}N\tau_n^2\|\alpha\|_{l^1(p)}^2+
\frac{1}{1-\nu}\left(\|P_{V_{\cB}}[\eps]\|_{n}^2+c_0^2\frac Sn \right)+\|\alpha\|_{l^1(\#(\cB))}^2\right)\\
&\leq \frac{4}{1-\nu}\;\#(\cI)\;\tau_n^2\,\|\alpha\|_{l^1(p)}^2
\left(\frac{\nu(1+\nu)}{1-\nu}\tau_n+1\right) + \frac{4\nu
c_0^2}{(1-\nu)^2}\,\tau_n \frac Sn + \frac{4
}{(1-\nu)^2}\;\#(\cI)\;\tau_n^2\|P_{V_{\cB}}[\eps]\|_{n}^2.
\end{align*}
This ends the proof.

%%%%%%%%%%%%%%%%%%%%%%%%%%%%%%%%%%%%%%%%%%%%%%%%%%%%%%%
\subsection{Proof of Proposition \protect{\ref{chi2bis}}}
%%%%%%%%%%%%%%%%%%%%%%%%%%%%%%%%%%%%%%%%%%%%%%%%%%%%%%%
First, we prove the part concerning the non random set $\cI$.
%%%%%%%%%%%%%%%%%%%%%%%%%%%%%%%%%%%%%%%%%%%%%%%%%%%%%%%
The following proposition gives the concentration inequalities when the
errors $\eps_i$'s are gaussian. Note that a corresponding inequality
stating concentration for projections of subgaussian variables can
be found in Proposition 5.1 (with possibly not the optimal constants
as stated by the authors) in \cite{huang-2009}.
\begin{lemma}\label{chi2}
Let $k$ be a positive integer and $U$ be a $\chi_k^2$ variable. Then
 $$\forall u^2\ge 4\frac kn,\quad P( \frac 1n\,U\ge u^2)\le \exp \left(-nu^2/8\right).$$
\end{lemma}
 Recall the following result by Massart (2007). If $X_t$ is be a centered gaussian process such that $\sigma^2:=\sup_t \E X_t^2$, then
\begin{equation}
\forall y>0,\quad P\left(\sup_t X_t-\E \sup_t X_t\ge y\right)\le
\exp -\frac {y^2}{2\sigma^2}.\label{supgaus}
\end{equation}
Let $Z_1,\ldots,Z_k$ i.i.d. standard Gaussian variables such that
\begin{align*}
P( U\ge nu^2)&=P(\sum_{i=1}^kZ_i^2\ge n u^2)
= P(\sup_{a\in S_1} \sum_{i=1}^ka_iZ_i\ge (nu^2)^{1/2})\\
&= P\left(\sup_{a\in S_1} \sum_{i=1}^ka_iZ_i-\E\sup_{a\in S_1}
\sum_{i=1}^ka_iZ_i\ge (nu^2)^{1/2}-\E\sup_{a\in S_1}
\sum_{i=1}^ka_iZ_i\right)
\end{align*}
where $S_1=\{a\in R^k, \|a_i\|_{l^2(k)}=1\}$. Denote
$$X_a=\sum_{i=1}^ka_iZ_i\quad\mbox{ and }\quad y=(nu^2)^{1/2}-\E\sup_{a\in S_1} \sum_{i=1}^ka_iZ_i.$$
Notice that
$$a\in S_1 \Rightarrow \E\left(X_a\right)^2=1$$  as well as
$$\E\sup_{a\in S_1}X_a=\E\left[\sum_{i=1}^kZ_i^2\right]^{1/2}\le \left[\E\sum_{i=1}^kZ_i^2\right]^{1/2}=k^{1/2}. $$
Since $u^2\ge 4\frac kn$, the announced result is proved as soon as
$y> (nu^2)^{1/2}/2$.

Assume now that $\cI$ is random and take into account all the non
random possibilities $\cI^\prime$  for the set
  $\cI$ and applying Proposition \ref{chi2bis} in the non random case. We get
 \begin{align*}
 P\left(\frac{1}{\sigma^2\,n}\|P_{V_{\cI}}[\eps]\|_{l^2(n)}^2\geq
 \mu^2\right)&\leq \sum_{\cI^\prime\subset\{1,\ldots,p\}}
 P\left(\frac{1}{\sigma^2\,n}\|P_{V_{\cI^\prime}}[\eps]\|_{l^2(n)}^2\geq
 \mu^2\right)\\
 &\leq p^{n_{\cI}}\exp\left(-n\mu^2/8\right)\\
&\leq
\exp\left(-n\mu^2\left[1/8-\frac{n_{\cI}\log{p}}{n\mu^2}\right]\right)\\
&\leq \exp\left(-n\mu^2/16\right)
 \end{align*}
as soon as $\mu^2\geq 16\;n_{\cI}\log{p}/n$.

%%%%%%%%%%%%%%%%%%%%%%%%%%%%%%%%%%%%%%%%%%%%%%%
\bibliographystyle{chicago}
\bibliography{BiblioLOL}

\vspace{1cm}

Address for correspondence:  [ PICARD Dominique, La\-bo\-ra\-toi\-re
de pro\-ba\-bi\-li\-t\'es et mo\-d\`e\-les al\'ea\-toi\-res, 175,
rue du Che\-va\-le\-ret, 75013 Paris, France].
picard@math.jussieu.fr

\newpage

%%%%%%%%%%%%%%%%%%%%%%%%%%%%%%%%%%%%%%%%%%%%%%%%%%%%%%%%%%%%%%%%%%%%%%%%%%%%%%%%%%%%%%%%%%%%%%%%%%%%%%%%%%%%%%%%%%%%
\begin{table}
\caption{Prediction error for varying sparsities $S$ and different
 distributions for the regressors, $n=250,p=1000$. $SNR=5$. \label{tablaws}}
\begin{tabular}{|c|c|c|c|c|c|c|c|c|}
  \hline
  S&  G & U & B & T(5) & T(4) & T(2) & T(1) \\
    \hline
5  & 0.00 (0.0)& 0.00 (0.00)& 0.00 (0.00)& 0.00 (0.00)&  0.00 (0.00)& 0.00 (0.00)& 0.00 (0.01) \\
\hline
10  & 0.00 (0.01)& 0.00 (0.02)& 0.00 (0.01)& 0.00 (0.00)& 0.00 (0.01)& 0.00 (0.00)& 0.00 (0.05) \\
\hline
15 & 0.01 (0.02)& 0.02 (0.03)& 0.02 (0.02)& 0.03 (0.03)& 0.01 (0.02)& 0.01 (0.02)& 0.01 (0.07) \\
\hline
20  & 0.04 (0.03)& 0.03 (0.03)& 0.03 (0.03)& 0.05 (0.04)& 0.03 (0.03)& 0.03 (0.03)& 0.04 (0.12) \\
\hline
25  & 0.07 (0.04)& 0.07 (0.05)& 0.06 (0.04)& 0.06 (0.03)& 0.07 (0.04)& 0.07 (0.04)& 0.08 (0.14) \\
\hline
30  & 0.10 (0.06)& 0.11 (0.06)& 0.08 (0.03)& 0.08 (0.04)& 0.11 (0.05)& 0.10 (0.05)& 0.17 (0.24) \\
\hline
35  & 0.15 (0.06)& 0.14 (0.07)& 0.15 (0.08)& 0.14 (0.06)& 0.13 (0.06)& 0.16 (0.07)& 0.25 (0.26) \\
\hline
40  & 0.19 (0.07)& 0.17 (0.06)& 0.17 (0.07)& 0.17 (0.06)& 0.18 (0.07)& 0.21 (0.09)& 0.35 (0.27) \\
  \hline
\end{tabular}
\end{table}

%%%%%%%%%%%%%%%%%%%%%%%%%%%%%%%%%%%%%%%%%%%%%%%%%%%%%%%%%%%%%%%%%%%%%%%%%
% Comparaison des mÈthodes
\begin{table}
\caption{Prediction errors for varying sparsity $S$ and varying SNR
computed using LOL, SIS-Reg and SIS-Lasso procedures.
 $n=200$, $p=1000$, $K=100$. \label{methodstudy}}
\begin{tabular}{|l|l|l|l|l|l|l|}
  \hline
  SNR & method  & $S=10$ & $S=20$ & $S=30$ & $S=50$ & $S=60$\\
  \hline

10 & LOL & 0.146 (0.141) &  0.273 (0.110) &  0.381 (0.068)  &  0.491 (0.118) &  0.462 (0.108) \\
10 & SIS-Lasso & 0.161 (0.103) &  0.389 (0.035) &  0.477 (0.030)  &  0.543 (0.029) &  0.554 (0.028) \\
10 & Lasso-Reg & 0.096 (0.005) &  0.095 (0.005) &  0.165 (0.102)  &  0.486 (0.121) &  0.472 (0.101) \\
\hline
5 & LOL & 0.228 (0.073) &  0.351 (0.077) &  0.436 (0.123)  &  0.478 (0.093) &  0.496 (0.067) \\
5 & SIS-Lasso & 0.223 (0.048) &  0.388 (0.053) &  0.476 (0.029)  &  0.543 (0.030) &  0.562 (0.032) \\
5 & Lasso-Reg & 0.188 (0.011) &  0.192 (0.016) &  0.323 (0.090)  &  0.466 (0.095) &  0.523 (0.124) \\
\hline
2 & LOL & 0.388 (0.071) &  0.463 (0.084) &  0.472 (0.060)  &  0.560 (0.150) &  0.545 (0.104)  \\
2 & SIS-Lasso & 0.418 (0.035) &  0.509 (0.026) &  0.541 (0.031)  &  0.589 (0.033) &  0.613 (0.032)  \\
2  & Lasso-Reg & 0.459 (0.052) &  0.514 (0.069) &  0.523 (0.065)  &  0.581 (0.153) &  0.597 (0.112)\\
\hline
\end{tabular}
\end{table}

%%%%%%%%%%%%%%%%%%%%%%%%%%%%%%%%%%%%%%%%%%%%%%%%%%%%%%%%%%%%%%%%%%%%%%%%%
% Very high dimension
\begin{table}
\caption{Prediction errors for varying ultra  high dimension $p$ and
sparsity $S$ computed using LOL. $SNR=5$,
$K=100$.\label{tabveryhigh}}
\begin{tabular}{|l|l|l|l|l|l|l|}
  \hline
  p & n & \multicolumn{5}{c|}{ S } \\
  \cline{3-7}
   & & 5 & 10 & 20 & 40 & 60\\
   \hline
     5000 & 400 & 0.195 (0.007) & 0.194 (0.006) & 0.236 (0.051) & 0.426 (0.058) & 0.497 (0.065)\\
      & 800 & 0.195 (0.004) & 0.195 (0.005) & 0.196 (0.012) & 0.234 (0.036) & 0.340 (0.046)\\
  \hline
  10000 & 400 & 0.195 (0.008) & 0.193 (0.007) & 0.244 (0.064) & 0.420 (0.052) & 0.443 (0.068)\\
      & 800 & 0.196 (0.004) & 0.195 (0.005) & 0.193 (0.005) & 0.236 (0.043) & 0.348 (0.050)\\
  \hline
  20000 & 400 & 0.204 (0.063) & 0.201 (0.049) & 0.277 (0.088) & 0.408 (0.074) & 0.401 (0.074)\\
      & 800 & 0.193  (0.004)& 0.195  (0.005)& 0.194  (0.004)& 0.242  (0.036) & 0.395 (0.055)\\
  \hline
\end{tabular}
\end{table}

\end{document}